\documentclass[11pt]{article}
\usepackage[a4paper,margin=26mm]{geometry}
\usepackage[T1]{fontenc}
\usepackage{lmodern,microtype}
\usepackage{amsmath,amssymb,amsthm,mathtools,mathrsfs,booktabs,enumitem}
\usepackage{xcolor}
\usepackage[colorlinks=true,linkcolor=blue,citecolor=blue,urlcolor=blue,hyperfootnotes=false]{hyperref}
\hypersetup{pdftitle={Asymptotically good codes from Siegel modular varieties},pdfauthor={Yota Maeda}}
\newif\ifshowrevisions
\showrevisionsfalse
\DeclareRobustCommand{\rev}[1]{\ifshowrevisions
  {\color{red}\hypersetup{allcolors=red}#1}\else#1\fi}
\newif\ifrevisionreference
\newcommand{\revisionbibstyle}[1]{%
  \ifshowrevisions
    \def\revisionkey{#1}%
    \def\revisionkeyA{AtiyahSchmid}%
    \def\revisionkeyB{LanCompactifications}%
    \revisionreferencefalse
    \ifx\revisionkey\revisionkeyA\revisionreferencetrue\fi
    \ifx\revisionkey\revisionkeyB\revisionreferencetrue\fi
    \ifrevisionreference
      \color{red}\hypersetup{allcolors=red}%
    \else
      \color{black}\hypersetup{linkcolor=blue,citecolor=blue,urlcolor=blue}%
    \fi
  \fi}
\makeatletter
\def\revision@lbibitem[#1]#2{%
  \revisionbibstyle{#2}%
  \revision@original@lbibitem[#1]{#2}}
\AtBeginDocument{%
  \let\revision@original@lbibitem\@lbibitem
  \let\@lbibitem\revision@lbibitem}
\makeatother
\newtheorem{theorem}{Theorem}[section]
\newtheorem{proposition}[theorem]{Proposition}
\newtheorem{lemma}[theorem]{Lemma}
\newtheorem{corollary}[theorem]{Corollary}
\theoremstyle{definition}
\theoremstyle{remark}\newtheorem{remark}[theorem]{Remark}
\newcommand{\F}{\mathbb F}
\newcommand{\Z}{\mathbb Z}
\newcommand{\Q}{\mathbb Q}
\newcommand{\R}{\mathbb R}
\newcommand{\A}{\mathcal A}
\newcommand{\PP}{\mathcal P}
\newcommand{\CC}{\mathscr C}
\newcommand{\YY}{\mathscr Y}
\newcommand{\Ha}{\operatorname{Ha}}
\newcommand{\ev}{\operatorname{ev}}
\newcommand{\wt}{\operatorname{wt}}
\newcommand{\Sp}{\operatorname{Sp}}
\newcommand{\Div}{\operatorname{div}}
\newcommand{\im}{\operatorname{im}}
\newcommand{\I}{\mathcal I}
\title{Asymptotically good codes from Siegel modular varieties}
\author{Yota Maeda\thanks{E-mail: \href{mailto:yota.maeda@mosk.tytlabs.co.jp}{\texttt{yota.maeda@mosk.tytlabs.co.jp}}}\\[3pt]
\small Toyota Central R\&D Labs., Inc., Japan}
\date{}

\begin{document}
\maketitle
\begin{abstract}
We construct asymptotically good codes by evaluating Siegel modular forms
at superspecial points on congruence towers of Siegel modular varieties.
For every genus $g\ge2$ and every fixed weight $w\ge g+1$, the construction
gives such a family in sufficiently large characteristic, over a finite
field which is independent of the level of the tower.
The towers can be chosen so that this field is $\F_{p^2}$.
We give an explicit lower bound for the asymptotic rate.
The lower bound for the relative minimum distance follows from the
Pieri formula on the space of flags in the Hodge bundle and an
intersection argument along its strata.
We also prove that, for $p\ge5$ and $0\le w\le p+1$, the kernel
of evaluation in every genus consists precisely of the multiples of the
ordinary Hasse invariant.
For $g=2,3$, the geometry of the supersingular locus gives
explicit bounds for the relative minimum distance and recovery sets
whose sizes are bounded independently of the level.
\end{abstract}
\section{Introduction}
\begingroup
\renewcommand\thefootnote{}
\footnotetext{\textit{2020 Mathematics Subject Classification.}
94B27, 11F46, 14G35.\newline
\textit{Keywords.} Siegel modular varieties, algebraic geometry codes,
superspecial points, Hasse invariant, Hecke orbits.}
\endgroup

Algebraic geometry codes are obtained by evaluating functions, or
sections of a line bundle, at rational points of an algebraic variety
over a finite field. On a curve, the degree of the line bundle bounds
the number of zeros of a nonzero section, while the Riemann--Roch
theorem controls the dimension of the space of sections.
Tsfasman, Vl\u{a}du\c{t}, and Zink used towers of modular or Shimura
curves to construct asymptotically good codes \cite{TVZ}.
It is natural to ask whether congruence towers of modular varieties
of higher dimension admit a similar construction.

In this paper, we consider Siegel modular varieties of arbitrary
genus $g\ge2$. We evaluate scalar Siegel modular forms at their
superspecial points and prove that, for a fixed weight and sufficiently
large characteristic, the resulting codes are asymptotically good as
the level increases. The finite field is fixed throughout the tower.
The main issue is to control the vanishing of modular forms at the
evaluation points uniformly in the level. We do this in two steps:
we determine the kernel of evaluation in an explicit range of weights,
and then bound the number of nonzero coordinates by means of the
stratification of the space of Hodge flags.

\subsection{Main results}

Let $p\ge5$ and $\ell\ne p$ be primes, and let $N\ge3$ be prime to
$p\ell$. Put
\[
 a_N^\pm=\min\{a\ge1:p^a\equiv\pm1\pmod N\},
 \qquad k=\F_{p^{2a_N^\pm}}.
\]
Fix a primitive $N$th root of unity in $k$.
For $r\ge1$, let
$X_{g,r}=\A_{g,N,\Gamma_0(\ell^r)}$
be the fine moduli scheme of principally polarized abelian varieties
of dimension $g$ with full symplectic level $N$ and a Lagrangian
direct summand of the $\ell^r$-torsion. The latter is an isotropic
subgroup which is, \'etale locally, a direct summand of rank $g$
over $\Z/\ell^r\Z$.
Let $\PP_{g,r}$ be its set of geometric superspecial points and put
$n_{g,r}=|\PP_{g,r}|$. These points are all rational over $k$, at
every level $r$; see Proposition~\ref{prop:field}.
\rev{The terms superspecial, supersingular, and ordinary are recalled
in Section~\ref{sec:geometry}.}
In particular, $k=\F_{p^2}$ is sufficient if $N\mid p-1$ or
$N\mid p+1$. Thus the fixed choice $N=3$ gives $k=\F_{p^2}$
for every $p\ge5$.

Choose a smooth projective toroidal compactification
$\overline X_{g,r}$, let $\mathbb E$ be its extended Hodge bundle,
whose fiber at $A$ is $H^0(A,\Omega_A^1)$, and put
$\omega=\det\mathbb E$. For an integer $w\ge0$, define
\[
 M_{g,r,w}=H^0(\overline X_{g,r},\omega^w),\qquad
 C_{g,r,w}=\im\bigl(M_{g,r,w}\xrightarrow{\ev}k^{\PP_{g,r}}\bigr).
\]
Here $\omega^w$ denotes $\omega^{\otimes w}$. Evaluation uses a
choice of basis in each fiber; changing the bases rescales the
coordinates and does not change the dimension or the minimum
distance. We put $M_{g,r,w}=0$ for $w<0$.

For a word $c$, its \textit{Hamming weight} $\wt(c)$ is the number of nonzero
coordinates. A nonzero code $C\subset k^n$ has \textit{minimum distance}
$d(C)=\min_{0\ne c\in C}\wt(c)$, \textit{rate} $\dim_k C/n$, and \textit{relative
minimum distance} $d(C)/n$. A family over one finite field is
\emph{asymptotically good} if its lengths tend to infinity and its
rates and relative minimum distances have positive lower limits.
In our families, $g,p,N,\ell$, and $w$ are fixed, and only $r$ varies.
Write
\[
 Q_g(p)=\prod_{j=1}^g\bigl(p^j+(-1)^j\bigr),\qquad
 R_{g,w}(p)=\frac1{Q_g(p)}
 \prod_{1\le i\le j\le g}\frac{2w-i-j}{2g+2-i-j}.
\]

\begin{theorem}\label{thm:general}
Fix $g\ge2$, $N\ge3$, and an integer $w\ge g+1$.
There exists $p_0(g,N,w)$ such that, for every prime
$p>p_0(g,N,w)$ with $p\nmid N$ and every prime $\ell\nmid pN$,
with $\ell$ odd if $g\ge3$,
the codes $C_{g,r,w}$ are asymptotically good over $k$ as
$r\to\infty$. Evaluation is injective at every level, and
\begin{equation}\label{eq:rate-general-weight}
 \liminf_{r\to\infty}\frac{\dim_k C_{g,r,w}}{n_{g,r}}
 \ge R_{g,w}(p)>0.
\end{equation}
There is a constant $\delta_{g,N,p,w}>0$, independent of $r$ and
$\ell$, such that
\[
 \frac{d(C_{g,r,w})}{n_{g,r}}\ge\delta_{g,N,p,w}>0
\]
whenever $C_{g,r,w}\ne0$.
\end{theorem}

One may take
\begin{equation}\label{eq:characteristic-threshold}
 p_0(g,N,w)=
 \begin{cases}
  w+1,&g=2,\\
  \max\{w+1,K+Cw\},&g\ge3,
 \end{cases}
\end{equation}
\rev{where $K$ and $C$ depend only on $g$ and $N$ and are chosen
in Lemma~\ref{lem:flag-boundary}.}
The distance constant $\delta_{g,N,p,w}$ is given by
\eqref{eq:distance-two} for $g=2$ and
\eqref{eq:distance-general} for $g\ge3$.
At the smallest weight and the next weight, the rate bounds are
\begin{equation}\label{eq:rate-general-main}
 R_{g,g+1}(p)=\frac1{Q_g(p)},\qquad
 R_{g,g+2}(p)=\frac{\binom{2g}{g}-\binom{2g}{g-2}}{Q_g(p)}.
\end{equation}
For general $g$, the constants in the distance bound depend on
the auxiliary ample class constructed in Section~\ref{sec:distance};
we do not give numerical values for these constants.
In genera two and three we obtain more explicit estimates, valid
also in small characteristics; these are given in
Proposition~\ref{prop:explicit}.

Our second main result describes the kernel of evaluation.
Let $S_{g,r}\subset\overline X_{g,r}$ be the reduced supersingular
locus. It is proper and lies in $X_{g,r}$.
The ordinary Hasse invariant
$\Ha\in H^0(\overline X_{g,r},\omega^{p-1})$
vanishes precisely on the closure of the nonordinary locus.
Thus every multiple of $\Ha$ vanishes at the superspecial points.
In small weights the converse holds in every genus.

\begin{theorem}\label{thm:mainideal}
Let $g\ge2$, $p\ge5$, and $0\le w\le p+1$.
At every level $r\ge1$,
\begin{equation}\label{eq:kernel-general}
 \ker(\ev)
 =H^0(\overline X_{g,r},\I_{S_{g,r}}\otimes\omega^w)
 =\Ha M_{g,r,w-(p-1)}.
\end{equation}
In particular, evaluation is injective for $w<p-1$ and has
one-dimensional kernel $k\Ha$ for $w=p-1$.
\end{theorem}

Theorem~\ref{thm:mainideal} holds for every $p\ge5$, independently
of $g$, and allows $\ell=2$ in every genus.
It describes the kernel on the whole space of scalar modular
forms in the stated range. In particular,
\[
 \dim_k C_{g,r,w}
 =\dim_k M_{g,r,w}-\dim_k M_{g,r,w-(p-1)}.
\]
For $g=2,3$, we also prove that vanishing on the entire supersingular
locus implies divisibility by $\Ha$ for
$0\le w<(p^2-1)/2$; see Proposition~\ref{prop:ideal2} and
Proposition~\ref{prop:ideal3}.

\subsection{Related work}

The construction of codes on varieties of higher dimension has been
studied. Hansen's method of auxiliary curves
\cite{Hansen} bounds zeros by restricting sections to curves through
the evaluation points. One must also control the curves on which
a section vanishes identically. Couvreur, Lebacque, and Perret
explain that this control is not automatically preserved under
finite covers \cite[Sections~1.7.2 and~1.9.2]{CLP}.
They give a different distance bound preserved under pullback
by finite morphisms and a criterion for infinite towers of \'etale covers
in which prescribed points split completely. They do not,
however, obtain new asymptotically good code families from
that criterion; finding towers and divisors that yield both
positive asymptotic rate and positive relative distance
remains an additional requirement
\cite[Introduction, Section~1.10, and Theorem~2.7]{CLP}.

\rev{Products of curves give tensor product codes
\cite[Proposition~5.1]{BlacheHallouin}, and hence yield
asymptotically good families from those on curves. Blache and Hallouin also
construct such families on ruled surfaces, with asymptotic
parameter bounds improving those of the corresponding product
codes \cite[Section~5]{BlacheHallouin}.}
Their construction uses ruled surfaces over a sequence of
curves. The present paper instead studies congruence towers
of Siegel modular varieties, using the kernel of evaluation
and the stratification of the Hodge flag space.

The superspecial mass formula is given in \cite{YuMass}.
In characteristic zero, the results of Finis--Lapid on limit
multiplicities \cite{FinisLapidApprox,FinisLapid} and the dimension
formula of Wakatsuki \cite{Wakatsuki} give the required asymptotic
dimension formula. Reduction modulo $p$ then gives a lower bound
for the dimensions in characteristic $p$.

Restriction to superspecial points was used by Ghitza
\cite{Ghitza} to compare systems of Hecke eigenvalues.
Yamauchi \cite{Yamauchi} studied restriction to strata and
divisibility by partial Hasse invariants in genus two.
The geometric argument underlying the case $g=2$ of
Theorem~\ref{thm:mainideal} already occurs in the proof of
\cite[Theorem~8.10]{CalegariGeraghty}.
\rev{We extend this argument to every genus by restricting to products
and using the density of Hecke orbits in central leaves.}
The density of prime-to-$p$ Hecke orbits in central leaves is
proved for Shimura varieties of Hodge type in
\cite[Theorem~I and Theorem~8.3.2]{DAddezioVanHoften}.

Ekedahl and van der Geer \cite{EkedahlGeer} study the
stratification of the Hodge flag space and establish normality,
irreducibility criteria, and a Pieri formula for its strata.
Local recovery by interpolation on curves is studied in
\cite{BargTamoVladut,BargEtAl}.

\subsection{Outline of the proofs}

For Theorem~\ref{thm:mainideal}, let $f$ vanish at all superspecial
points. Consider the largest closed subset of its zero locus
which is stable under prime-to-$p$ Hecke correspondences.
It contains the superspecial locus.
Restriction to a product of an abelian surface and fixed elliptic
curves allows us to apply the calculation in genus two.
The density of Hecke orbits then propagates vanishing through
central leaves with successively larger ordinary parts.
The last of these is the locus of $p$-rank $g-1$.
Its closure is the reduced divisor of $\Ha$, so $f$ is divisible
by $\Ha$.

For the distance estimate in genus $g\ge3$, we use an integral
closed flag stratum dominating the locus of $p$-rank zero.
A form with nonzero evaluation is nonzero on this stratum.
The Pieri formula gives a lower bound, linear in $p$, for the
ratio of the sum of the degrees of its boundary divisors to
the degree of the stratum.
We apply the same estimate successively to strata of smaller
dimension, ending at flags over superspecial points.
All auxiliary classes are chosen at the fixed level $N$ and
pulled back through the tower. This makes the resulting relative
bound independent of $r$ and $\ell$.
Theorem~\ref{thm:mainideal} supplies injectivity for $w<p-1$;
the estimate for the dimension of $M_{g,r,w}$ then proves
Theorem~\ref{thm:general}.

\subsection{Organization of the paper}

Section~\ref{sec:geometry} treats the levels, the field of
definition of the evaluation points, and the Hasse divisor.
Section~\ref{sec:kernel} proves Theorem~\ref{thm:mainideal}.
Section~\ref{sec:distance} establishes the estimates for the
minimum distance, and Section~\ref{sec:growth} proves
Theorem~\ref{thm:general}.
Section~\ref{sec:examples} gives explicit bounds and local
recovery in genera two and three.
The calculations for the supersingular locus in genus three
are collected in Appendix~\ref{sec:genus3}.

\medskip
\noindent\textbf{Notation.}
We write $\lambda_i=c_i(\mathbb E)$ and
$\lambda=\lambda_1=c_1(\omega)$.
For a $k$-scheme, sheaf, or section, a subscript $\overline{k}$
denotes extension of scalars to an algebraic closure of $k$.
The spaces of modular forms and the codes are defined over $k$;
geometric components and intersection products are considered
over $\overline{k}$.
In particular,
\[
 M_{g,r,w}\otimes_k\overline{k}
 \simeq H^0((\overline X_{g,r})_{\overline{k}},
                 \omega_{\overline{k}}^w).
\]
All loci defined by conditions on abelian varieties carry their
reduced structures. For a closed subscheme $Y\subset X$,
we denote its ideal sheaf by $\I_Y$.
\section{Levels, compactifications, and the Hasse divisor}\label{sec:geometry}

An abelian variety is \textit{supersingular} if, over an algebraic closure,
it is isogenous to a product of supersingular elliptic curves.
It is \textit{superspecial} if it is isomorphic to such a product.
These isogenies and isomorphisms are not required to preserve
the polarization. The \textit{$p$-rank} of an abelian variety $A$ of
dimension $g$ is the integer $f$ such that
$|A[p](\overline{k})|=p^f$; $A$ is \textit{ordinary} if $f=g$.
At a fixed level, the superspecial locus is finite, whereas the
supersingular locus has positive dimension for $g\ge2$.

\subsection{The congruence tower}

Let $\Sp_{2g}$ be the symplectic group for the standard alternating form.
Write $\Gamma_g(N)$ for the kernel of the reduction map
$\Sp_{2g}(\Z)\to\Sp_{2g}(\Z/N\Z)$,
and let $\Gamma_{0,g}(\ell^r)$ be the subgroup whose lower-left $g\times g$ block is zero modulo $\ell^r$. Put $\Gamma_{g,r}=\Gamma_g(N)\cap\Gamma_{0,g}(\ell^r)$. Let $d_g=g(g+1)/2=\dim X_{g,r}$.
For an integer $M\ge3$ with $p\nmid M$, let $\A_{g,M}$
denote the Siegel moduli scheme over $\overline{k}$
with full symplectic level $M$.
We choose the primitive roots of unity defining the
symplectic pairings compatibly.
This moduli scheme corresponds to the principal
congruence subgroup $\Gamma_g(M)$.
The inclusions
$\Gamma_g(N\ell^r)\subset\Gamma_{g,r}\subset\Gamma_g(N)$
induce forgetful morphisms
\[
 \A_{g,N\ell^r}
 \to (X_{g,r})_{\overline{k}}=X_{g,r}\times_k\overline{k}
 \to \A_{g,N}.
\]
We first determine a common field of definition for the superspecial points of $X_{g,r}$.
\begin{proposition}\label{prop:field}
With $a_N^\pm$ and $k$ as in the introduction, all geometric
superspecial points of $X_{g,r}$ are rational over
$k=\F_{p^{2a_N^\pm}}$ for every $r\ge1$.
In particular, one may take $k=\F_{p^2}$ if $N\mid p-1$ or
$N\mid p+1$, and hence whenever $N=3$.
\end{proposition}
\begin{proof}
By \cite[Proposition~5.1 and Definition~5.2]{YuArithmetic}, every
principally polarized superspecial abelian variety admits a model
over $\F_{p^2}$ whose Frobenius endomorphism is $[-p]$.
Extend this model to $k$.
Its Frobenius endomorphism is $[c]$, where $c=(-p)^{a_N^\pm}$ and
$c\equiv\pm1\pmod N$.

If $c\equiv1\pmod N$, retain this model. Otherwise, take its
quadratic twist by the automorphism $[-1]$ over $k$.
This automorphism preserves the polarization, so the polarization
descends to the twist and remains principal. After identifying
the twist with the original variety over $\overline{k}$, its
Frobenius endomorphism is $[-c]$. Thus in either case we have
a model whose Frobenius is scalar and acts trivially on the
$N$-torsion. Every full level $N$ structure is therefore defined
over $k$. The chosen pairing is also defined over $k$, since
$p^{2a_N^\pm}\equiv1\pmod N$.

On the $\ell^r$-torsion, Frobenius acts by an invertible scalar.
It therefore preserves every Lagrangian direct summand, which
descends to a subgroup scheme over $k$. The resulting object
with both level structures defines a $k$-point of the fine
moduli scheme. This argument is independent of $r$.
\end{proof}

Note that the field in Proposition~\ref{prop:field} is a sufficient field of definition; it need not be minimal. 
We next compute the code lengths.
Set
\begin{equation}\label{eq:index}
 D_{g,r}=|\Sp_{2g}(\Z/N\Z)|\,
 \ell^{d_g(r-1)}\prod_{j=1}^{g}(\ell^j+1).
\end{equation}
This is the index of $\Gamma_g(N)\cap\Gamma_{0,g}(\ell^r)$ in $\Sp_{2g}(\Z)$. Indeed the Lagrangian Grassmannian over $\F_\ell$ has $\prod_{j=1}^{g}(\ell^j+1)$ points. Its smoothness and dimension $d_g$ show that every Lagrangian direct summand at level $\ell^r$ has $\ell^{d_g}$ lifts at the next level. The factors at $N$ and $\ell$ are independent. In particular $D_{g,r}\to\infty$.

By Proposition~\ref{prop:field}, all geometric superspecial
points are $k$-rational.
Since $N\ge 3$, a principally polarized
abelian variety with full symplectic level $N$ structure
has no nontrivial automorphisms.
Thus, the superspecial mass formulas
\cite[Theorem~3.1 and Lemma~3.2]{YuMass}, together with
\eqref{eq:index}, give, in every genus,
\begin{equation}\label{eq:mass-general}
 n_{g,r}=D_{g,r}v_gQ_g(p),\qquad
 v_g=2^{-g}\prod_{j=1}^{g}|\zeta(1-2j)|.
\end{equation}
Here $\zeta$ is the Riemann zeta function and
$v_2=1/5760$, $v_3=1/2903040$. In particular,
\begin{equation}\label{eq:mass}
 n_{2,r}=D_{2,r}\frac{(p-1)(p^2+1)}{5760},\qquad
 n_{3,r}=D_{3,r}
 \frac{(p-1)(p^2+1)(p^3-1)}{2903040}.
\end{equation}
In particular, $n_{g,r}\to\infty$, and the ratio
$n_{g,r}/D_{g,r}$ is independent of $r$.

\subsection{Compactification and the Hasse divisor}\label{subsec:hasse}

\rev{Fix a discrete valuation ring $R$ of mixed characteristic with
residue field $k$, containing a primitive $N$th root of unity that
lifts the chosen root in $k$.
For each $r$, choose a smooth projective toroidal compactification
$\overline{\mathcal X}_{g,r}$ over $R$, with extended Hodge bundle,
and take $\overline X_{g,r}$ to be its special fiber.
These compactifications and extensions are provided by
\cite[Chapters IV and V]{FC}; see also
\cite[Theorem~6.4.1.1, Proposition~7.3.1.4, and
Theorem~7.3.3.4]{LanCompactifications}.
Write $D_{\mathrm{tor}}$ for the boundary divisor on the compactification
and for its restrictions to the fibers.}
The line bundle $\omega$ is nef and big: a positive tensor power is the pullback of an ample line bundle on the minimal compactification
\rev{\cite[Theorem~7.2.4.1(2),(3)]{LanCompactifications}}.
In particular $\lambda^{d_g}>0$. All intersection identities below are taken on a toroidal compactification, with rational coefficients, and may be computed after extending the field to $\overline{k}$. We identify top-dimensional intersection products with their degrees.

For the universal abelian scheme $A$, Verschiebung $V:A^{(p)}\to A$ induces $V^*: \mathbb E\to \mathbb E^{(p)}$, where $(p)$ denotes pullback by the absolute Frobenius of the base. The \textit{ordinary Hasse invariant} $\Ha$ is its determinant, viewed as a section of $\omega^{p-1}$ via $\det(\mathbb E^{(p)})\simeq\omega^p$. This section extends to the toroidal compactification. We require its divisor with multiplicities.

Let $H\subset\overline X_{g,r}$ be the reduced closure, over $k$,
of the nonordinary locus of $X_{g,r}$.
Since $k$ is perfect, $H_{\overline{k}}$ is the reduced closure
of the nonordinary locus of $(X_{g,r})_{\overline{k}}$.
Thus $H_{\overline{k}}\cap(X_{g,r})_{\overline{k}}$ is the locus
where the $p$-rank is at most $g-1$, with its reduced structure.
\begin{lemma}\label{lem:hasse}
For every $g\ge2$, the divisor of $\Ha$ is reduced and geometrically integral. It is the closure of the nonordinary locus and satisfies
\begin{equation}\label{eq:hasse}
 \Div(\Ha)=H,\qquad [H]=(p-1)\lambda.
\end{equation}
\end{lemma}
\begin{proof}
\rev{At full prime-to-$p$ level, the locus of $p$-rank $g-1$
is dense in the nonordinary locus; see \cite[Introduction]{YuMono}.
It is irreducible by the argument in the proof of
\cite[Proposition~3.1]{YuMono}, using the irreducibility theorem
\cite[Theorem~3.2]{YuMono}; here $g-1>0$ excludes the
supersingular case.} Its image under the surjective
map $\A_{g,N\ell^r}\to(X_{g,r})_{\overline{k}}$ is therefore
irreducible. Thus $H_{\overline{k}}$ is integral.

On the open moduli scheme, the divisor of $\Ha$ is
$mH_{\overline{k}}$ for some $m\ge1$.
In genus two, the local calculation at a point of $p$-rank one
in \cite[Proposition~3.24]{Yamauchi} gives $m=1$.
For $g>2$, take a formal curve \rev{with parameter space
$\operatorname{Spec}\overline{k}[[t]]$}
through such a point on which
the Hasse invariant in genus two has order one, and take the product
of the corresponding family with a fixed ordinary abelian variety
of dimension $g-2$.
\rev{Since $N\ell$ is prime to $p$, the required torsion is finite
\'etale and hence constant over the strictly henselian ring
$\overline{k}[[t]]$. Compatible full symplectic levels therefore
extend over this curve and define a family in $X_{g,r}$ after
forgetting part of the level.}
The Hasse invariant of a product is the product of the Hasse
invariants, so its pullback to this formal curve again has order
one. Since the \rev{moduli scheme} is smooth, a local equation
for $\Ha$ is a unit times the $m$th power of an equation for
$H_{\overline{k}}$. The order along the curve forces $m=1$.

There is no additional component supported on the toroidal
boundary. Indeed, at the generic point of a boundary divisor,
the universal semi-abelian variety is an extension of a generic
abelian variety of smaller dimension by a torus
\cite[Chapter~IV]{FC}.
\rev{The boundary stratum dominates the corresponding moduli
space of abelian varieties of smaller dimension
\cite[Theorem~7.2.4.1(4),(5)]{LanCompactifications}. Its abelian
part is therefore ordinary at the generic point, since the
ordinary locus is dense in that space.}
Verschiebung induces an isomorphism on the invariant
differentials of the torus. Its determinant is therefore nonzero.

The zero scheme of $\Ha$ is a Cartier divisor on a smooth variety.
It is Cohen--Macaulay, has no embedded associated points, and is
generically reduced by the preceding calculation. It is thus
reduced, and its support is $H_{\overline{k}}$.
Faithful flatness gives $\Div(\Ha)=H$ over $k$, and
$[H]=c_1(\omega^{p-1})=(p-1)\lambda$.
\end{proof}

\rev{Since $\omega$ is nef and big, there are no nonzero sections of
$\omega^u$ for $u<0$: the divisor of such a section would have
negative intersection with $\lambda^{d_g-1}$.
Moreover, $M_{g,r,0}=k$ since $\overline X_{g,r}$ is proper and
geometrically integral.}

We also record the relation between restriction to $H$ and
divisibility by the Hasse invariant.
By Lemma~\ref{lem:hasse}, $H=\Div(\Ha)$ is a reduced effective
Cartier divisor, with ideal sheaf
$\I_H=\mathcal O_{\overline X_{g,r}}(-H)\simeq\omega^{-(p-1)}$.
The exact sequence for restriction to $H$, tensored with
$\omega^w$, therefore gives
\begin{equation}\label{eq:hasse-restriction}
 \ker\bigl(M_{g,r,w}\to H^0(H,\omega^w|_H)\bigr)
 =\Ha M_{g,r,w-(p-1)}.
\end{equation}
Thus a modular form vanishes on $H$ if and only if it is divisible
by $\Ha$.
Section~\ref{sec:kernel} shows that, in small weights, vanishing
at the superspecial points already forces vanishing on $H$.

\section{The kernel of evaluation}\label{sec:kernel}

We first recall the calculation in genus two. We then use Hecke
correspondences to deduce Theorem~\ref{thm:mainideal} in arbitrary
genus. We allow finer full level structures in this section,
as needed for the pullbacks in Section~\ref{subsec:hecke-products}.

\subsection{The calculation in genus two}\label{sec:genus2}

The geometric supersingular locus in genus two is a union of
projective lines. Every component $C$ contains $p^2+1$
superspecial points, exactly $p+1$ components pass through each
superspecial point, and
\begin{equation}\label{eq:g2degree}
 \deg(\omega_{\overline{k}}|_C)=p-1;
\end{equation}
see \cite[Theorem~1.1 and Section~4.4]{YuMass},
\cite[Section~2]{KO}, and \cite[Section~2, following Theorem~2.10]{vdG}.
\rev{These statements hold at our levels and at any finer full
prime-to-$p$ level: the inverse image of a component under a
finite \'etale map which forgets part of the level is a disjoint
union of projective lines, and the Hodge bundle pulls back
under this map. The calculation below therefore applies at
all these levels.}
Since the supersingular locus is the locus of $p$-rank zero,
its reduced cycle class is
\begin{equation}\label{eq:g2cycles}
 [(S_{2,r})_{\overline{k}}]
 =(p-1)(p^2-1)\lambda_{2,\overline{k}}
 =\frac{(p-1)(p^2-1)}2\lambda_{\overline{k}}^2.
\end{equation}
Here we used $\lambda_1^2=2\lambda_2$; see
\cite[Section~2]{vdGH}.

\begin{proposition}\label{prop:ideal2}
For $g=2$ and $0\le w<(p^2-1)/2$,
\[
 H^0(\overline X_{2,r},\I_{S_{2,r}}\otimes\omega^w)
 =\Ha M_{2,r,w-(p-1)}.
\]
For $0\le w\le p+1$, this space also equals $\ker(\ev)$.
\end{proposition}
\begin{proof}
The inclusion from right to left follows from $S_{2,r}\subset H$.
Conversely, let $f$ vanish on $S_{2,r}$ and suppose that
$f_{\overline{k}}|_{H_{\overline{k}}}\ne0$.
The zero divisor of this restriction contains every component of
$(S_{2,r})_{\overline{k}}$ with multiplicity at least one.
Intersecting with the nef class $\lambda_{\overline{k}}$ gives
\[
 \lambda_{\overline{k}}[(S_{2,r})_{\overline{k}}]
 \le w\lambda_{\overline{k}}^2[H_{\overline{k}}].
\]
By \eqref{eq:hasse} and \eqref{eq:g2cycles}, and since
$\lambda_{\overline{k}}^3>0$, this implies
$w\ge(p^2-1)/2$, a contradiction.
Thus $f|_H=0$, and \eqref{eq:hasse-restriction} proves the
first assertion.

Now suppose $w\le p+1$ and $f\in\ker(\ev)$.
If its restriction to a supersingular component $C$ were nonzero,
it would have $p^2+1$ distinct zeros but degree
$w(p-1)\le(p+1)(p-1)=p^2-1$.
Hence it vanishes on every component, and therefore on the reduced
supersingular locus. Since $p+1<(p^2-1)/2$ for $p\ge5$, the
first assertion applies.
\end{proof}

\subsection{Hecke correspondences and restriction to products}
\label{subsec:hecke-products}

\rev{Fix a Siegel modular form $f\in M_{g,r,w}$. Throughout the remainder of
this section, the geometric constructions are over $\overline{k}$.}
For background on Hecke correspondences, see
\cite[Sections~1.2--1.7]{ChaiOort}.
\rev{In the Siegel case, \textit{Hecke correspondences} parametrize isogenies
$\varphi:(A,\lambda_A)\to(A',\lambda_{A'})$ between principally
polarized abelian varieties with compatible level structures such that}
$\varphi^*\lambda_{A'}=n\lambda_A$ for a positive integer $n$.
We refer to prime-to-$p$ correspondences if $p\nmid n$.
The two projections record the source and target, and the \textit{Hecke orbit}
of a point consists of the points obtained by these correspondences.

\rev{Let $\operatorname{GSp}_{2g}$ be the group of similitudes
of the standard alternating form. Write $\mathbb A_f$ for the
ring of finite adeles of $\Q$ and $\mathbb A_f^p$ for the same
restricted product with the $p$-component omitted.}
To use all prime-to-$p$ Hecke correspondences, we
temporarily allow all primitive roots defining the symplectic
pairing at full level $N$. Thus we work on the full Siegel
moduli scheme $\mathscr X$ for the corresponding compact open
subgroup of $\operatorname{GSp}_{2g}(\mathbb A_f^p)$, with
\rev{hyperspecial level $\operatorname{GSp}_{2g}(\Z_p)$ at $p$}.
The scheme $(X_{g,r})_{\overline{k}}$ is one open and closed
component of $\mathscr X$.
Extend $f_{\overline{k}}$ by zero on the other components and
denote the resulting section again by $f$.
This convention avoids restricting the Hecke correspondences
to those preserving one choice of the pairing.

Let $\mathcal O^p(x)$ be the prime-to-$p$ Hecke orbit of a
geometric point $x\in\mathscr X$. Put
\begin{equation}\label{eq:hecke-zero}
 Z(f)=\{x\in\mathscr X:f(y)=0\text{ for every }
                         y\in\mathcal O^p(x)\}.
\end{equation}
We give this set its reduced structure.

\begin{lemma}\label{lem:hecke-zero}
The set $Z(f)$ is closed and stable under prime-to-$p$ Hecke
correspondences. It is contained in the zero locus of $f$.
If $f$ vanishes at all superspecial points, then $Z(f)$ contains
the superspecial locus.
\end{lemma}
\begin{proof}
Write a Hecke correspondence as
$\mathscr X\xleftarrow{q_1}T\xrightarrow{q_2}\mathscr X$.
Both projections are finite \'etale. If $D$ is the zero locus
of $f$, the condition that every point in this correspondence
above $x$ maps into $D$ defines the closed set
$\mathscr X\setminus q_1\bigl(q_2^{-1}(\mathscr X\setminus D)\bigr)$,
because $q_1$ is open. Intersecting over all correspondences,
including their compositions and the identity, gives $Z(f)$.
The description by entire orbits proves its stability and its
containment in $D$.
A prime-to-$p$ isogeny induces an isomorphism of $p$-divisible
groups and preserves superspeciality. This proves the last assertion.
\end{proof}

Let $B$ be a fixed principally polarized abelian variety of
dimension $g-2$. Choose a full level $M$, prime to $p$ and
sufficiently divisible that the product construction defines a morphism
$j:\A_{2,M}\to\mathscr X$, $A\mapsto A\times B$,
using the product polarization and chosen compatible levels.
\rev{Let $H_{2,M}$ denote the nonordinary locus in the open
moduli scheme $\A_{2,M}$.}

\begin{lemma}\label{lem:product-propagation}
Suppose $0\le w\le p+1$. If $j$ sends every superspecial point
of $\A_{2,M}$ into $Z(f)$, then $j(H_{2,M})\subset Z(f)$.
\end{lemma}
\begin{proof}
Fix a prime-to-$p$ Hecke correspondence
$\mathscr X\xleftarrow{q_1}T\xrightarrow{q_2}\mathscr X$
and pull it back along $j$.
After replacing $M$ by a sufficiently divisible integer $M'$ prime
to $p$, this pullback is a trivial finite \'etale cover:
\[
 \A_{2,M'}\times_{\mathscr X,q_1}T
 \simeq\coprod_\alpha\A_{2,M'}.
\]
Here the map from $\A_{2,M'}$ to $\mathscr X$ is the product map
with the additional level structure forgotten. Composing each
component with $q_2$ gives a morphism
$j_\alpha:\A_{2,M'}\to\mathscr X$.
Indeed, the isogenies and the target level structures are specified
by subgroups and bases in torsion of some order prime to $p$.
A full level of a sufficiently divisible order trivializes these
finite data. The torsion of the fixed factor $B$ is constant.
Thus every branch is defined after adding a full level structure.
The integer $M'$ may depend on the chosen Hecke correspondence.

The universal prime-to-$p$ isogeny identifies the Hodge bundles
on the source and target. The Hodge bundle of $A\times B$ is
the direct sum of those of its factors. In particular,
\[
 j_\alpha^*\omega_g\simeq
 \omega_2\otimes
 \bigl(\det H^0(B,\Omega_B^1)\otimes_{\overline{k}}\mathcal O\bigr),
\]
where $\omega_g$ is the Hodge line bundle on $\mathscr X$.
Choosing a basis of the fixed line $\det H^0(B,\Omega_B^1)$ identifies
$j_\alpha^*f$ with a section of $\omega_2^w$ on $\A_{2,M'}$.
The boundary of the minimal compactification of $\A_{2,M'}$
has codimension two. Koecher's principle therefore extends this
section to a modular form on a toroidal compactification;
see \cite[Theorem~2.3 and Remark~2.4]{LanKoecher}.

Every superspecial point of $\A_{2,M'}$ maps under $j$ into
$Z(f)$, and all its Hecke images are zeros of $f$.
Consequently $j_\alpha^*f$ vanishes at all superspecial points.
Proposition~\ref{prop:ideal2} implies that it vanishes on
$H_{2,M'}$.
This holds for every branch of every Hecke correspondence.
Since the map $\A_{2,M'}\to\A_{2,M}$ is surjective, all Hecke
images of $j(H_{2,M})$ are zeros of $f$, which is the required
inclusion.
\end{proof}

\subsection{Proof of Theorem~\ref{thm:mainideal}}
\label{subsec:kernel-proof}

For an abelian variety $A$, the \textit{$p$-divisible group}
$A[p^\infty]=\varinjlim_n A[p^n]$ records its group schemes of
$p$-power torsion. A principal polarization on $A$ induces an
isomorphism $A[p^\infty]\to A^\vee[p^\infty]$, called a \textit{principal
quasipolarization}; see \cite[Definition~10.6, Section~10.7,
and Corollary~2.14]{ChaiOort}.
\rev{An isomorphism of principally quasipolarized $p$-divisible
groups is required to respect these maps to the dual groups,
or equivalently the induced alternating Weil pairings.}
A \textit{central leaf} is the reduced locally closed locus on which
this isomorphism class is fixed; see
\cite[Section~1.1.1]{DAddezioVanHoften}.

Fix an ordinary elliptic curve $E_{\mathrm{ord}}$ and a
supersingular elliptic curve $E_{\mathrm{ss}}$ over $\overline{k}$.
Let $G_{\mathrm{ord}}$ and $G_{\mathrm{ss}}$ be their principally
quasipolarized $p$-divisible groups.
For $0\le a\le g-1$, let $\mathcal C_a\subset\mathscr X$ be the
central leaf defined by $G_{\mathrm{ord}}^a\oplus G_{\mathrm{ss}}^{g-a}$.
Thus it contains the moduli points of abelian varieties whose
$p$-divisible groups, with their quasipolarizations, are isomorphic
to that of $E_{\mathrm{ord}}^a\times E_{\mathrm{ss}}^{g-a}$.
The product groups above admit similitudes with every multiplier
in $\Z_p^\times$. For $G_{\mathrm{ord}}$, this follows from its
decomposition into multiplicative and \'etale parts. For
$G_{\mathrm{ss}}$, it follows from the surjectivity of the reduced
norm on the units of its quaternion endomorphism ring; see
\cite[Corollaries~2.6 and~2.9]{Ghitza} for the description of this
ring and its action on the polarization. Hence multiplying the
polarization by a $p$-adic unit does not change these isomorphism
classes, so the loci $\mathcal C_a$ are stable under prime-to-$p$
Hecke correspondences.

We use the density of prime-to-$p$ Hecke orbits in these leaves,
as proved in \cite[Theorem~I and Theorem~8.3.2]{DAddezioVanHoften}.
\rev{The Newton slopes on $\mathcal C_a$ are $0$ with multiplicity
$a$, $1/2$ with multiplicity $2g-2a$, and $1$ with multiplicity
$a$. In the Siegel case, \textit{basic} means that all slopes are $1/2$.
Thus $\mathcal C_a$ is nonbasic for $a\ge1$.
The adjoint group of $\operatorname{GSp}_{2g}$ is $\Q$-simple,
so the hypotheses of that theorem hold.} In particular, a nonempty closed subset
of $\mathcal C_a$ stable under these correspondences is all of
$\mathcal C_a$.

The leaf $\mathcal C_0$ is contained in the superspecial locus.
Moreover, every principally polarized superspecial abelian
surface has the same principally quasipolarized $p$-divisible
group as $G_{\mathrm{ss}}^2$; see \cite[Corollary~2.8]{Ghitza}.
Finally, $\mathcal C_{g-1}$ is the locus of $p$-rank $g-1$:
the \'etale and multiplicative parts split over an algebraically
closed field, and \rev{the remaining principally quasipolarized
group is connected and has connected dual. It has height two
and dimension one, so it is isomorphic to $G_{\mathrm{ss}}$}.

\begin{proof}[Proof of Theorem~\ref{thm:mainideal}]
The assertion for $g=2$ is Proposition~\ref{prop:ideal2}.
Let $g\ge3$ and let $f\in\ker(\ev)$.
Use the extension to $\mathscr X$ described above.
Lemma~\ref{lem:hecke-zero} gives
$\mathcal C_0\subset Z(f)$.
Suppose inductively that $\mathcal C_a\subset Z(f)$ for some
$0\le a\le g-2$, and take
$B=E_{\mathrm{ord}}^a\times E_{\mathrm{ss}}^{g-a-2}$.
For the product map $j:A\mapsto A\times B$, the superspecial
points of $\A_{2,M}$ map into $\mathcal C_a$.
Lemma~\ref{lem:product-propagation} therefore gives
$j(H_{2,M})\subset Z(f)$.
Taking $A=E_{\mathrm{ord}}\times E_{\mathrm{ss}}$ shows that
$Z(f)\cap\mathcal C_{a+1}$ is nonempty.
This intersection is closed in the central leaf and stable
under prime-to-$p$ Hecke correspondences.
The density theorem gives
$\mathcal C_{a+1}\subset Z(f)$.

By induction, $f$ vanishes on $\mathcal C_{g-1}$, and hence on
the nonordinary divisor of the original component.
Since this divisor is reduced, the vanishing is scheme-theoretic.
It descends to $k$ by faithful flatness.
Equation~\eqref{eq:hasse-restriction} now gives
$f\in\Ha M_{g,r,w-(p-1)}$.
The reverse inclusions
\[
 \Ha M_{g,r,w-(p-1)}
 \subset H^0(\overline X_{g,r},\I_{S_{g,r}}\otimes\omega^w)
 \subset\ker(\ev)
\]
are immediate. This proves all equalities.
The last assertions follow from the absence of modular forms
of negative weight and from $M_{g,r,0}=k$.
\end{proof}
\section{Estimates for the minimum distance}\label{sec:distance}

We begin with two intersection estimates, then apply them to the
Hodge flag space for $g\ge3$. The estimate in genus two is given
at the end of the section. All varieties and intersection products
in this section are over $\overline{k}$, and we suppress this
subscript on the Hodge classes. We continue to evaluate forms
defined over $k$ at the original $k$-rational points.

\subsection{Intersection and incidence}

Let $X$ be a projective scheme, let $L$ be a nef line bundle,
and let \rev{$w$ be a nonnegative integer}.
For an integral subvariety $Y$ of dimension $t$,
write $\deg_L(Y)=\deg(c_1(L)^t\cap[Y])$.

\begin{lemma}\label{lem:intersection}
Assume $\dim Y=t\ge1$, and let $s\in H^0(Y,L^w)$ be nonzero. If $W_1,\ldots,W_h$ are distinct integral divisors on $Y$ and $s$ vanishes on each $W_i$, then
\[
 \sum_{i=1}^{h}\deg_L(W_i)\le w\deg_L(Y).
\]
\end{lemma}
\begin{proof}
The zero divisor $D$ of $s$ is an effective Cartier divisor
with $\mathcal O_Y(D)\simeq L^w|_Y$. As cycles,
$[D]-\sum_i[W_i]$ is effective, since each $W_i$ occurs in
$D$ with multiplicity at least one. Intersecting with
$c_1(L|_Y)^{t-1}$ and using nefness gives
\[
 \sum_i\deg_L(W_i)
 \le\deg\bigl(c_1(L|_Y)^{t-1}\cap[D]\bigr)
 =w\deg_L(Y).
\]
\end{proof}

To iterate this estimate, it is useful to measure degrees by an
auxiliary ample class. The following proposition does not require
the strata of a given dimension to have equal degrees.

\begin{proposition}[Incidence and degrees]
\label{prop:weighted-incidence}
Let $X$ be a projective scheme over an algebraically closed field,
let $L$ be a nef line bundle, and let $A$ be an ample integral
divisor class such that $A-c_1(L)$ is nef.
For an integral $t$-dimensional subvariety $Y$, write
$\deg_A(Y)=\deg(A^t\cap[Y])$.
For an integer $d\ge1$ and $0\le t\le d$, let $\mathcal T_t$ be a finite collection of
integral $t$-dimensional subvarieties, with
$\mathcal T_d=\{Z\}$ and $\mathcal T_0$ a set of points.
Let $w\ge0$ be an integer. For each $1\le t\le d$, suppose
there are constants $a_t>0$ and $b_t>w$ such that, for every
$Y\in\mathcal T_t$,
\begin{equation}\label{eq:weighted-boundary}
 \sum_{\substack{W\in\mathcal T_{t-1}\\W\subset Y}}
       \deg_A(W)\ge b_t\deg_A(Y),
\end{equation}
and each $W\in\mathcal T_{t-1}$ is contained in at most $a_t$
members of $\mathcal T_t$.
If $s\in H^0(X,L^w)$ and $s|_Z\ne0$, then
\begin{equation}\label{eq:weighted-incidence}
 |\{P\in\mathcal T_0:s(P)\ne0\}|
 \ge \deg_A(Z)\prod_{t=1}^d\frac{b_t-w}{a_t}>0.
\end{equation}
In particular, if $a_t=a$ and $b_t=b$ for every $t$, this becomes
\[
 |\{P\in\mathcal T_0:s(P)\ne0\}|
 \ge\left(\frac{b-w}{a}\right)^d\deg_A(Z).
\]
\end{proposition}
\begin{proof}
For $Y\in\mathcal T_t$ with $s|_Y\ne0$, its zero divisor gives
\[
 \sum_{\substack{W\in\mathcal T_{t-1},\ W\subset Y\\s|_W=0}}
       \deg_A(W)
 \le w\deg\bigl(c_1(L)A^{t-1}\cap[Y]\bigr)
 \le w\deg_A(Y).
\]
The first inequality follows from the argument with effective
divisors in Lemma~\ref{lem:intersection}; the second uses the nefness of
$A-c_1(L)$ and $A$.
Consequently the sum of the degrees of the divisors on which
$s$ is nonzero is at least $(b_t-w)\deg_A(Y)$.
Put
\[
 S_t=\sum_{\substack{Y\in\mathcal T_t\\s|_Y\ne0}}\deg_A(Y).
\]
Summing the preceding inequality and using the bound on the
number of containing strata gives $a_tS_{t-1}\ge(b_t-w)S_t$.
Since $S_d=\deg_A(Z)$ and a geometric point has $A$-degree one,
iteration proves \eqref{eq:weighted-incidence}.
\end{proof}

\subsection{Strata on the Hodge flag space}\label{sec:higher}

Assume $g\ge3$.
Let $Z_{g,r}$ be the reduced locus of $p$-rank zero in
$(X_{g,r})_{\overline{k}}$. It has dimension $g(g-1)/2$ and is proper:
\rev{in the extension of the stratification to a toroidal
compactification, the \textit{extended $p$-rank} at a boundary point is
the sum of the torus rank and the $p$-rank of the abelian part.
It is therefore positive at every boundary point; see
\cite[Section~5 and Lemma~5.1]{EkedahlGeer}.}
The Hodge class is ample on $Z_{g,r}$: it comes from an ample
class on the minimal compactification, and the map to the minimal
compactification is an isomorphism over $X_{g,r}$.

Let $\pi_r:\mathcal F_r\to(X_{g,r})_{\overline{k}}$ be the
bundle of complete flags in the Hodge bundle $\mathbb E$.
\rev{Let $a:\mathcal A_r\to(X_{g,r})_{\overline{k}}$
be the universal abelian scheme.
Its de Rham bundle
$\mathcal H_r=R^1a_*\Omega^\bullet_{\mathcal A_r/(X_{g,r})_{\overline{k}}}$
has a symplectic pairing induced by the principal polarization,
and $\mathbb E\subset\mathcal H_r$ is a Lagrangian subbundle.
Verschiebung induces a surjective map
$V:\mathcal H_r\to\mathbb E^{(p)}$, where $(p)$ denotes
pullback by absolute Frobenius.
On $\mathcal F_r$, the tautological flag
$0=E_0\subset E_1\subset\cdots\subset E_g=\pi_r^*\mathbb E$
extends to a symplectic flag by $E_{g+i}=E_{g-i}^{\perp}$.
Its conjugate flag is defined by
\[
 D_g=\ker V,\qquad D_{g+i}=V^{-1}(E_i^{(p)}),\qquad
 D_{g-i}=D_{g+i}^{\perp}\quad(0\le i\le g),
\]
where the bundles and $V$ are pulled back to $\mathcal F_r$.
We identify the Weyl group $W(C_g)$ with the permutations $u$
of $\{1,\ldots,2g\}$ satisfying
$u(2g+1-i)=2g+1-u(i)$.
The stratum $U_{u,r}$ is defined by the conditions
\[
 \dim(E_i\cap D_j)=\#\{a\le i:u(a)\le j\}
 \quad(0\le i,j\le2g)
\]
on geometric fibers; see \cite[Sections~3--4]{EkedahlGeer}.}
Write $\operatorname{len}(u)$ for the length of $u$ and
$v\prec u$ for a Bruhat cover, meaning
$v<u$ and $\operatorname{len}(v)=\operatorname{len}(u)-1$.
By \cite[Corollary~8.4]{EkedahlGeer},
$\dim U_{u,r}=\operatorname{len}(u)$ and
$\overline U_{u,r}=\bigcup_{v\le u}U_{v,r}$,
and every closed stratum $\overline U_{u,r}$ is reduced,
normal, and Cohen--Macaulay. In particular, its irreducible
components are disjoint.
These statements are preserved by the finite \'etale maps
that forget part of the level structure. They therefore hold
at our levels by descent from full level $N\ell^r$.

\rev{The subgroup $S_g\subset W(C_g)$ stabilizes
$\{1,\ldots,g\}$. The minimal-length representatives of the
cosets $uS_g$ are called \textit{final elements}. For such an element,
its final type is determined by
\[
 \nu_u(i)=i-\#\{a\le g:u(a)\le i\},\qquad 1\le i\le g.
\]
The associated Young diagram has row lengths given by the
positive values of
$\#\{1\le i\le g:i-\nu_u(i)\ge j\}$, for $1\le j\le g$;
we use the conventions of \cite[Section~2.2]{EkedahlGeer}.
Let $u_0$ be the final element with $\nu_{u_0}(i)=i-1$.
Its diagram is the single row $[g]$. Set
$\widetilde Z_r=\overline U_{u_0,r}$.
The Ekedahl--Oort stratification classifies principally polarized
abelian varieties by the isomorphism type of their polarized
$p$-torsion group schemes.
The image of $U_{u_0,r}$ is the dense Ekedahl--Oort stratum
of $Z_{g,r}$, and the map to that image is finite, surjective,
and \'etale.} At full level $N\ell^r$, 
\cite[Theorem~11.5]{EkedahlGeer} implies that the stratum of type $u_0$ is
connected and hence irreducible, since it is smooth.
Its image under the surjective morphism that forgets part of
the level structure is $U_{u_0,r}$, which is therefore irreducible.
The criterion applies because, for $g\ge3$, the
diagram $[g]$, which consists of one row, omits a row of some length between
$\lceil(g+1)/2\rceil$ and $g$.
Thus $\widetilde Z_r$ is integral of dimension $g(g-1)/2$, and
$\pi_r:\widetilde Z_r\to Z_{g,r}$ is proper and dominant.
It is projective, since it is closed in the flag bundle over
the proper variety $Z_{g,r}$.

For $0\le t\le g(g-1)/2$, let $\mathcal T_{t,r}$ consist of all
irreducible components of $\overline U_{u,r}$ with $u\le u_0$
and $\operatorname{len}(u)=t$.
Then $\mathcal T_{g(g-1)/2,r}=\{\widetilde Z_r\}$.
If $e$ is the identity of $W(C_g)$, the points of
$\mathcal T_{0,r}=U_{e,r}$ lie over the superspecial points.
By \cite[Lemma~4.6 and Corollary~8.4(iii)]{EkedahlGeer},
there is an integer $h_g(p)>0$, independent of the point and
the prime-to-$p$ level, such that every superspecial point has
exactly $h_g(p)$ geometric points in $U_{e,r}$ above it.
Consequently
\begin{equation}\label{eq:flag-evaluations}
 |\mathcal T_{0,r}|=h_g(p)n_{g,r},\qquad
 |\{Q\in\mathcal T_{0,r}:(\pi_r^*f)(Q)\ne0\}|
 =h_g(p)\wt(\ev(f)).
\end{equation}

\subsection{Degrees of the boundary strata}

Here the \textit{boundary of a stratum} means its boundary in its closure
in the Hodge flag space, rather than the toroidal boundary.
We choose the auxiliary ample class at the fixed full level $N$.
All classes in what follows are pulled back from that level to
$\mathcal F_r$ and then restricted to its proper closed strata.
Let $E_\bullet$ be the tautological Hodge flag, put
$x_i=c_1(E_i/E_{i-1})$, and write $h=\pi_r^*\lambda$ and $L_i=-\sum_{j=1}^{i}x_j\quad(1\le i\le g)$ with $L_0=0$.
Thus $L_g=-h$; the $L_i$ are the fundamental weights in the
convention of \cite[Section~10]{EkedahlGeer}.
The Weyl group acts on the $x_i$ by signed permutations.
For an expression such as $uL_i$, this action is taken on the
weight lattice before passing to the corresponding divisor class.

\cite[Proposition~10.3(i)--(ii)]{EkedahlGeer} gives integers
$b_i\ge1$, chosen independently of $p$ at the fixed level $N$,
such that the line bundles with classes
$N_i=L_i+b_i h$, $1\le i<g$, are generated by global sections.
Put
\begin{equation}\label{eq:auxiliary-classes}
 \begin{aligned}
 b&=\sum_{i=1}^{g-1}b_i,\\
 A&=\sum_{i=1}^{g-1}N_i+(g-1)h
   =\sum_{i=1}^{g-1}L_i+(b+g-1)h,\\
 B&=\sum_{i=1}^{g-1}L_i+bL_g.
 \end{aligned}
\end{equation}
The class $A$ is ample on the flag space and hence on its
proper closed strata. Moreover, $A-h$ is nef there.
Pullback under the finite maps that forget part of the level
structure preserves these properties. Thus the same choices
work for every $r$ and every admissible $\ell$.

\begin{lemma}\label{lem:flag-boundary}
There are constants $K\ge0$ and $C>0$, depending only on $g$
and $N$, such that, for every $1\le t\le g(g-1)/2$ and
$Y\in\mathcal T_{t,r}$,
\begin{equation}\label{eq:flag-boundary}
 \sum_{\substack{W\in\mathcal T_{t-1,r}\\W\subset Y}}
       \deg_A(W)
 \ge\frac{p-K}{C}\deg_A(Y).
\end{equation}
Each $W\in\mathcal T_{t-1,r}$ is contained in at most $g^2$
members of $\mathcal T_{t,r}$.
\end{lemma}
\begin{proof}
The Pieri formula on the flag space
\cite[Theorem~10.1]{EkedahlGeer} is the identity
\begin{equation}\label{eq:flag-pieri}
 (ph+pL_{g-i}-uL_i)[\overline U_{u,r}]
 =\sum_{v\prec u}c^i_{u,v}[\overline U_{v,r}]
\end{equation}
in $\operatorname{CH}^1(\overline U_{u,r})_{\Q}$.
The coefficients $c^i_{u,v}$ are nonnegative, depend only on
the root system and the elements $u,v$, and are independent of $p$ and the level.
For each cover $v\prec u$, at least one $c^i_{u,v}$ is positive.
The independence of $p$ follows also directly from the
Chevalley coefficients in the proof of that theorem.
Summing \eqref{eq:flag-pieri} over $1\le i<g$ and adding $b$
times the formula for $i=g$ gives
\begin{equation}\label{eq:combined-pieri}
 (pA-uB)[\overline U_{u,r}]
 =\sum_{v\prec u}c_{u,v}[\overline U_{v,r}],\qquad
 c_{u,v}=\sum_{i=1}^{g-1}c^i_{u,v}+b c^g_{u,v}>0.
\end{equation}
Choose $C\ge1$ to bound these finitely many coefficients from above.

The classes $N_1,\ldots,N_{g-1},h$ form a basis of the rational
weight space, and all their coefficients in $A$ are positive.
Since the Weyl group is finite, one can choose $K\ge0$ so that
$KA-uB$ is a nonnegative rational linear combination of these
classes for every $u\in W(C_g)$.
It is therefore nef on every proper stratum.
The choices of $C$ and $K$ depend on $g,N$ and the fixed $b_i$,
but not on $p,r$, or $\ell$.

Normality makes each irreducible component $Y$ of
$\overline U_{u,r}$ open and closed in that closed stratum.
We may consequently restrict \eqref{eq:combined-pieri} to $Y$.
Its boundary terms are precisely the components $W$ of
$\overline U_{v,r}$ contained in $Y$, with the same coefficients
$c_{u,v}$.
Intersecting with $A^{t-1}$ gives
\[
 C\sum_{\substack{W\in\mathcal T_{t-1,r}\\W\subset Y}}
       \deg_A(W)
 \ge\deg\bigl((pA-uB)A^{t-1}\cap[Y]\bigr)
 \ge(p-K)\deg_A(Y),
\]
which proves \eqref{eq:flag-boundary}.

Finally, a fixed component of type $v$ lies in at most one
component of any given type $u$, because the components of
$\overline U_{u,r}$ are disjoint.
A cover $v\prec u$ is obtained using a reflection corresponding
to a positive root. The root system $C_g$ has $g^2$ positive
roots, so there are at most $g^2$ possible types of containing strata.
This proves the incidence bound.
\end{proof}

\subsection{The estimate in genus at least three}

Write $\widetilde Z_N$ for the integral flag stratum at full level
$N$, and let $n_{g,N}$ be the number of superspecial points at
that level. Thus
$n_{g,N}=|\Sp_{2g}(\Z/N\Z)|v_gQ_g(p)$.
\begin{proposition}\label{prop:distance-general}
Let $g\ge3$ and $w\ge0$. With the constants above, suppose
$p>K+Cw$. Every nonzero word $c\in C_{g,r,w}$ satisfies
\begin{equation}\label{eq:distance-general}
 \frac{\wt(c)}{n_{g,r}}
 \ge\delta_{g,N,p,w}
 :=\left(\frac{p-K-Cw}{Cg^2}\right)^{g(g-1)/2}
   \frac{\deg_A(\widetilde Z_N)}{h_g(p)n_{g,N}}>0.
\end{equation}
The bound is independent of $r$ and $\ell$.
\end{proposition}
\begin{proof}
Choose $f\in M_{g,r,w}$ with $\ev(f)=c\ne0$.
Since every superspecial point lies in $Z_{g,r}$, the restriction
of $f$ to $Z_{g,r}$ is nonzero. The proper dominant morphism
$\pi_r:\widetilde Z_r\to Z_{g,r}$ then gives a nonzero section
$s=\pi_r^*f|_{\widetilde Z_r}$.
\rev{Apply Proposition~\ref{prop:weighted-incidence} on the
projective variety $\widetilde Z_r$, with
$L=(\pi_r^*\omega)|_{\widetilde Z_r}$ and the collections
$\mathcal T_{t,r}$.}
The class $A-h$ is nef, and Lemma~\ref{lem:flag-boundary}
provides the constants \rev{$a_t=g^2$ and $b_t=(p-K)/C>w$
for every $t$}.
Consequently,
\[
 |\{Q\in\mathcal T_{0,r}:s(Q)\ne0\}|
 \ge\left(\frac{p-K-Cw}{Cg^2}\right)^{g(g-1)/2}
       \deg_A(\widetilde Z_r).
\]
By \eqref{eq:flag-evaluations}, the left-hand side is
$h_g(p)\wt(c)$.

The forgetful map to full level $N$ is finite \'etale of degree
$e_{g,r}=D_{g,r}/|\Sp_{2g}(\Z/N\Z)|$.
The flag strata and the class $A$ are pulled back under this map.
The projection formula and the mass formula therefore give
$\deg_A(\widetilde Z_r)=e_{g,r}\deg_A(\widetilde Z_N)$ and
$n_{g,r}=e_{g,r}n_{g,N}$.
Dividing by $h_g(p)n_{g,r}$ proves the assertion.
\end{proof}

\subsection{The estimate in genus two}

\begin{proposition}\label{prop:distance2}
For $g=2$ and $0\le w\le p+1$, every nonzero word
$c\in C_{2,r,w}$ satisfies
\begin{equation}\label{eq:distance-two}
 \frac{\wt(c)}{n_{2,r}}
 \ge\delta_{2,N,p,w}
 :=\left(1-\frac{2w}{p^2-1}\right)
       \left(1-\frac{w(p-1)}{p^2+1}\right)>0.
\end{equation}
\end{proposition}
\begin{proof}
Let $\mathscr C$ be the set of geometric supersingular components,
and choose $f$ with $\ev(f)=c\ne0$.
Then $f|_H\ne0$. By Lemma~\ref{lem:intersection}, the number
of components on which $f$ vanishes is at most
$w\lambda^2[H]/(p-1)$.
On the other hand, \eqref{eq:g2cycles} gives
\[
 |\mathscr C|(p-1)=\lambda[S_{2,r}]
 =\frac{(p-1)(p^2-1)}2\lambda^3.
\]
Thus the proportion of components on which $f$ is nonzero is
at least $1-2w/(p^2-1)$.
On each such component, at least $p^2+1-w(p-1)$ superspecial
points have nonzero evaluation.
Each superspecial point lies on $p+1$ components, and
$|\mathscr C|(p^2+1)=n_{2,r}(p+1)$.
Counting incidences gives the stated product.
Both factors are positive in the indicated range.
\end{proof}
\section{Dimensions of modular forms and asymptotic goodness}\label{sec:growth}

We now estimate the dimensions of spaces of modular forms as
the level increases. Together with Theorem~\ref{thm:mainideal}
and the distance estimates of Section~\ref{sec:distance}, this
will prove Theorem~\ref{thm:general}.

For a congruence subgroup
$\Gamma\subset\Sp_{2g}(\Z)$,
let $S_w(\Gamma)$ denote the space of scalar
holomorphic Siegel cusp forms of weight $w$.
We use asymptotic dimension formulas in characteristic
zero to obtain lower bounds for the corresponding
spaces of modular forms in characteristic $p$.

\begin{theorem}\label{thm:growth}
Fix $g\ge2$, $p,N,\ell$, and an integer $w\ge g+1$.
Assume that $\ell$ is odd if $g\ge3$. Then
\[
 \liminf_{r\to\infty}
 \frac{\dim_k M_{g,r,w}}{D_{g,r}}
 \ge\mu_{g,w},
\]
where
\begin{equation}\label{eq:mu-general}
 \mu_{g,w}
 =v_g\prod_{1\le i\le j\le g}
       \frac{2w-i-j}{2g+2-i-j}>0.
\end{equation}
In particular,
\begin{align}
 \mu_{2,w}&=\frac{(w-1)(w-2)(2w-3)}{34560},\label{eq:mu2}\\
 \mu_{3,w}&=\frac{(2w-2)(2w-3)(2w-4)^2(2w-5)(2w-6)}
 {2^{16}3^6 5^2 7}.\label{eq:mu3}
\end{align}
\end{theorem}
\begin{proof}
Put $G=\Sp_{2g}$ and $\widehat{\Z}=\prod_q\Z_q$, with $q$ ranging over primes. Let $K_r\subset G(\widehat{\Z})$ be the compact open subgroup imposing principal level $N$ at primes dividing $N$ and the stabilizer of the standard Lagrangian modulo $\ell^r$ at $\ell$. At all other primes its local factor is $G(\Z_q)$. Then
$G(\Q)\cap K_r=\Gamma_{g,r}$ and
$[G(\widehat{\Z}):K_r]=D_{g,r}$,
where the intersection is taken inside $G(\mathbb A_f)$. Normalize the finite Haar measure by $\operatorname{vol}(G(\widehat{\Z}))=1$ and fix a Haar measure at infinity.

For the $\Q$-simple simply connected group $G$, nondegeneracy in \cite[Definition~1.3]{FinisLapidApprox} reduces to $\operatorname{vol}(K_r)\to0$. Here $\operatorname{vol}(K_r)=D_{g,r}^{-1}\to0$, so this condition holds. Moreover $G(\Q)\cap K_r$ contains no nontrivial central element as $N\ge3$. The limit multiplicity theorem \cite[Corollary 2]{FinisLapid} therefore applies. For $\Sp_4$, use part~(1), concerning split groups of rank two. For $\Sp_{2g}$ with $g\ge3$, use part~(2), concerning quasi-split classical groups. The required bound on the level at $2$ holds because $\ell$ is odd and $N$ is fixed.

Let $\pi_w$ be the holomorphic discrete series of $G(\R)$ with
scalar lowest $\mathrm{U}(g)$-type $\det^w$. It exists for every $w>g$,
including $w=g+1$; see \cite[Sections~5.1 and~5.3]{Wakatsuki}.
\rev{Here $\mathrm{U}(g)$ is the standard maximal compact subgroup of $G(\R)$.}
Write $m(\pi_w,\Gamma_{g,r})$ for its multiplicity in the discrete
automorphic spectrum and $d(\pi_w)$ for its \textit{formal degree} with
respect to the chosen real Haar measure.
\rev{Thus $d(\pi_w)$ is the Plancherel mass of the singleton $\{\pi_w\}$.
Put $v_G=\operatorname{vol}(G(\Z)\backslash G(\R))$, so that
\[
 \operatorname{vol}(\Gamma_{g,r}\backslash G(\R))=D_{g,r}v_G.
\]
Since a discrete series representation is isolated in the tempered dual,
limit multiplicity applies to $\{\pi_w\}$. Together with strong
approximation, it gives}
\[
 \frac{m(\pi_w,\Gamma_{g,r})}{D_{g,r}}
 \to v_Gd(\pi_w)\rev{.}
\]
The realization of cusp forms as holomorphic discrete series gives
$m(\pi_w,\Gamma_{g,r})=\dim_{\mathbb C}S_w(\Gamma_{g,r})$
for $w>g$; see the beginning of \cite[Section~5.4]{Wakatsuki}.

\rev{The general formula of Harish--Chandra for the formal degree gives}
\[
 v_Gd(\pi_w)=c_g\prod_{1\le i\le j\le g}(2w-i-j)
 \qquad(w>g),
\]
where $c_g>0$ is independent of $w$;
\rev{see \cite[Theorem~8.1]{AtiyahSchmid}.
The same product appears in \cite[(5.9)]{Wakatsuki}, where the
discussion is restricted to the integrable discrete series. Here we
use the general formula, which applies throughout $w>g$.}
Indeed, the Harish-Chandra parameter is $(w-1,\ldots,w-g)$.
The factors from compact roots are independent of $w$, and the
remaining factors give the displayed product. All these parameters
lie in the same chamber when $w>g$.

To determine $c_g$, fix one integer $w_0\ge g+2$ and apply
the same limit multiplicity statement to principal levels $Nq^j$,
where $q$ is any fixed odd prime not dividing $N$.
At this weight, the dimension formula in
\cite[Theorem~1.6]{Wakatsuki} applies. After division by
$[G(\Z):\Gamma_g(Nq^j)]$, its limit as $j\to\infty$ is
\[
 2^{-3g}\prod_{j=1}^{g}
 \frac{(j-1)!\,|B_{2j}|}{(2j-1)!\,j}
 \prod_{1\le i\le j\le g}(2w_0-i-j),
\]
where $B_{2j}$ are the Bernoulli numbers.
Using $|\zeta(1-2j)|=|B_{2j}|/(2j)$ and
\[
 \prod_{1\le i\le j\le g}(2g+2-i-j)
 =2^g\prod_{j=1}^{g}\frac{(2j-1)!}{(j-1)!},
\]
we obtain $v_Gd(\pi_{w_0})=\mu_{g,w_0}$.
The formula for the formal degree then gives
$v_Gd(\pi_w)=\mu_{g,w}$ as in \eqref{eq:mu-general}
for every $w\ge g+1$.
Thus the dimension formula at finite level is used only at $w_0$,
within its stated range. Consequently
$\dim_{\mathbb C}S_w(\Gamma_{g,r})=D_{g,r}\mu_{g,w}+o(D_{g,r})$.

\rev{At each level use the smooth projective toroidal integral model
$\overline{\mathcal X}_{g,r}$ over $R$ fixed in
Section~\ref{subsec:hasse}, together with its boundary $D_{\mathrm{tor}}$.
We also write $\omega$ for the determinant of its extended Hodge bundle.
\textit{Cusp forms} are sections of $\omega^w(-D_{\mathrm{tor}})$.}
Its global sections form a finite torsion-free $R$-module whose
rank is the dimension of the space of cusp forms in characteristic
zero. If $\varpi$ is a uniformizer, the exact sequence for
multiplication by $\varpi$ gives an injection
\rev{\[
 H^0(\overline{\mathcal X}_{g,r},\omega^w(-D_{\mathrm{tor}}))\otimes_R k
 \to H^0(\overline X_{g,r},\omega^w(-D_{\mathrm{tor}}))
 \to M_{g,r,w}.
\]}
This proves the asserted lower bound. 
\end{proof}

\begin{proof}[Proof of Theorem~\ref{thm:general}]
Fix $g,N$, and $w$ as in the theorem.
Take $p_0(g,N,w)$ as in \eqref{eq:characteristic-threshold}.
Then $p>p_0(g,N,w)$ implies $p>w+1$ and, for $g\ge3$,
$p>K+Cw$.
Theorem~\ref{thm:mainideal} gives injectivity of evaluation,
since $w<p-1$. Theorem~\ref{thm:growth} and
\eqref{eq:mass-general} give
\[
 \liminf_{r\to\infty}\frac{\dim_k C_{g,r,w}}{n_{g,r}}
 \ge\frac{\mu_{g,w}}{v_gQ_g(p)}=R_{g,w}(p)>0.
\]
For $g\ge3$, Proposition~\ref{prop:distance-general} gives a
positive lower bound for the relative minimum distance,
independent of $r$ and $\ell$.
For $g=2$, use Proposition~\ref{prop:distance2}.
Since $n_{g,r}\to\infty$ and $k$ is fixed, the families are
asymptotically good.

The first equality in \eqref{eq:rate-general-main} follows by
substituting $w=g+1$. The second follows from
\[
 \prod_{1\le i\le j\le g}
 \frac{2g+4-i-j}{2g+2-i-j}
 =\frac1{g+2}\binom{2g+2}{g+1}
 =\binom{2g}{g}-\binom{2g}{g-2}.
\]
\end{proof}
\begin{remark}\label{rem:other-levels}
\rev{The same arguments apply if the Lagrangian direct summand is
replaced by an isotropic flag of locally free direct summands of
any fixed nonempty type $0<d_1<\cdots<d_s\le g$ in the
$\ell^r$-torsion. More precisely, let $P\subset\Sp_{2g}$ be the
parabolic subgroup scheme over $\Z$ stabilizing the standard
isotropic flag of this type, and replace}
$\Gamma_{0,g}(\ell^r)$ by the inverse image of
$P(\Z/\ell^r\Z)$ under reduction modulo $\ell^r$.
Keep the full level $N$ structure, and define the spaces of forms
and the evaluation codes as before. With $D_{g,r}$ replaced by
the corresponding index, Theorems~\ref{thm:general},
\ref{thm:mainideal}, and~\ref{thm:growth} hold with the same
rate bound and the same characteristic threshold.

To see this, scalar Frobenius preserves every such flag, so
Proposition~\ref{prop:field} gives the same common field $k$.
The maps that forget the flag are finite \'etale, and full level
$N\ell^r$ dominates these moduli schemes. The arguments for the
kernel, the irreducibility of the required flag stratum, and the
distance estimate therefore apply unchanged. Both the superspecial
count and the relevant degrees are multiplied by the degree of
the level map. Finally, if $d_P=\dim(\Sp_{2g}/P)>0$, the index is
\[
 |\Sp_{2g}(\Z/N\Z)|\,
 \ell^{d_P(r-1)}|(\Sp_{2g}/P)(\F_\ell)|,
\]
which tends to infinity. The limit multiplicity argument applies
with the same condition on the level at $2$, namely that $\ell$
be odd when $g\ge3$.
\end{remark}

\section{Explicit bounds and local recovery}\label{sec:examples}

The geometry of the supersingular locus gives explicit estimates
in genera two and three, including small characteristics.
The calculation for genus three is given in
Appendix~\ref{sec:genus3}.

\subsection{Explicit bounds}

\begin{proposition}
\label{prop:explicit}
The following families are asymptotically good over the fixed
finite field $k=\F_{p^{2a_N^\pm}}$ as $r\to\infty$.
\begin{enumerate}[label=\textup{(\roman*)},itemsep=5pt]
\item Suppose $g=2$, $p\ge5$, and $3\le w\le p-1$. Then
\begin{align}
 \liminf_{r\to\infty}
 \frac{\dim_k C_{2,r,w}}{n_{2,r}}
 &\ge
 \frac{(w-1)(w-2)(2w-3)}
 {6(p-1)(p^2+1)},
 \label{eq:rate2}\\
 \frac{d(C_{2,r,w})}{n_{2,r}}
 &\ge
 \left(1-\frac{2w}{p^2-1}\right)
 \left(1-\frac{w(p-1)}{p^2+1}\right).
 \label{eq:distance2}
\end{align}
\item Suppose $g=3$, $p\ge11$, $\ell$ is odd, and
$4\le w\le(p-1)/2$. Then evaluation is injective, and
\begin{align}
 \liminf_{r\to\infty}
 \frac{\dim_k C_{3,r,w}}{n_{3,r}}
 &\ge
 \frac{(2w-2)(2w-3)(2w-4)^2(2w-5)(2w-6)}
 {2880(p-1)(p^2+1)(p^3-1)},
 \label{eq:rate3}\\
 \frac{d(C_{3,r,w})}{n_{3,r}}
 &\ge
 \left(1-\frac{w}{(p-1)(p^2+1)}\right)
 \frac{(p^2+1-w(p-1))(p^2-p-1)}
 {(p^2+1)(p^2-p-1)+2w(p-1)}.
 \label{eq:distance3}
\end{align}
\end{enumerate}
The distance estimates hold at every level at which the code is
nonzero, including all sufficiently large $r$.
\end{proposition}

\begin{proof}
In genus two and $w<p-1$, Theorem~\ref{thm:mainideal} gives
$\dim C_{2,r,w}=\dim M_{2,r,w}$.
At the endpoint $w=p-1$, it gives
$\dim C_{2,r,p-1}=\dim M_{2,r,p-1}-1$, since $M_{2,r,0}=k$.
This constant loss does not affect the limiting rate.
Dividing the lower bound of Theorem~\ref{thm:growth}
by the mass formula in genus two, \eqref{eq:mass}, gives
\eqref{eq:rate2}.
The distance bound follows from Proposition~\ref{prop:distance2}.

In genus three, Proposition~\ref{prop:distance3} gives both
injectivity and the distance bound.
Dividing the bound of Theorem~\ref{thm:growth} by the
mass formula in genus three, \eqref{eq:mass}, gives
\eqref{eq:rate3}.
\rev{All stated lower bounds are positive, and
$n_{g,r}\to\infty$. Both families are therefore asymptotically good.}
\end{proof}

Table~\ref{tab:examples} lists the lower bounds from
Proposition~\ref{prop:explicit}. In every row, take $N=3$ and
a prime $\ell\nmid3p$, with $\ell$ odd in genus three.
The family is then defined over the fixed field $\F_{p^2}$.

\begin{table}[ht]
\centering
\begin{tabular}{ccccc}
\toprule
$g$ & $p$ & $w$ & Lower bound for the rate & Lower bound for $d/n$\\
\midrule
2&5&3&$1/104$&$21/52$\\
2&5&4&$5/104$&$10/39$\\
2&7&4&$1/60$&$13/30$\\
2&7&6&$1/10$&$21/100$\\
3&11&4&$1/1622600$&$1358576/2040145$\\
3&11&5&$1/115900$&$79461/136213$\\
3&13&5&$7/2239920$&$693935/1079976$\\
\bottomrule
\end{tabular}
\caption{Lower bounds for the asymptotic rate and the relative minimum distance.}
\label{tab:examples}
\end{table}

For a concrete tower, take $(g,p,N,\ell)=(2,5,3,7)$.
Then $k=\F_{25}$ and
\[
 n_{2,r}=374400\cdot7^{3(r-1)}.
\]
\rev{The first two rows of Table~\ref{tab:examples} give the bounds
at weights $3$ and $4$.}
For $p\to\infty$ and $w=\lfloor cp\rfloor$, the pairs of lower bounds
for the rate and relative minimum distance in
\eqref{eq:rate2}--\eqref{eq:distance2} and
\eqref{eq:rate3}--\eqref{eq:distance3} tend to
\[
 \begin{cases}
 (c^3/3,\ 1-c),&g=2,\quad 0<c<1,\\
 (c^6/45,\ 1-c),&g=3,\quad 0<c<1/2.
 \end{cases}
\]
These limits describe how the proved lower bounds depend on the ratio $w/p$. In genus three, counting vanishing orders along the contracted sections removes the factor $1-2c$ that would remain from the elementary incidence estimate alone. The powers $c^3$ and $c^6$ come from the weight dependence of the dimension estimates.
\rev{The field varies in this outer limit. Asymptotic goodness over
a fixed field refers to the level limit $r\to\infty$ with
$p,N,\ell$, and $w$ fixed.}

\subsection{Local recovery}\label{subsec:local-recovery}

A linear code has \textit{locality} at most $\rho$ if every coordinate is
a linear function of at most $\rho$ other coordinates on the code.
A \textit{recovery set} for a coordinate is a set of other coordinates
from which it can be recovered in this way.
More generally, a code has \textit{$(\rho,\delta_{\mathrm{loc}})$ locality}
if every coordinate belongs to a set of at most
$\rho+\delta_{\mathrm{loc}}-1$ coordinates such that every nonzero
word of the code restricted to this set has weight at least
$\delta_{\mathrm{loc}}$.
Local recovery by interpolation on curves is used in
\cite{BargTamoVladut,BargEtAl}; here it applies to the curves already
used in the distance estimates.
\rev{For $g=2$, let $\CC$ be the irreducible components of the
geometric supersingular locus, as in Section~\ref{sec:genus2}.
For $g=3$, let $\CC$ be the curves in
Proposition~\ref{prop:geometry3}. In both cases, every evaluation
point lies on such a curve, and each curve is a projective line
with $\deg(\omega_{\overline{k}}|_C)=p-1$ containing
$p^2+1$ superspecial points.}
For $C\in\CC$, put $\PP_C=C\cap\PP_{g,r}$, using the geometric
images of the $k$-rational evaluation points, and write
$C_{g,r,w}|_{\PP_C}$ for the code
obtained by retaining only these coordinates.

\begin{corollary}\label{cor:locality}
Let $g=2,3$, $0\le w\le p+1$, and $t=w(p-1)$.
Each set $\PP_C$ has $p^2+1$ elements, and
\[
 \dim_k(C_{g,r,w}|_{\PP_C})\le t+1.
\]
Every nonzero word of this restricted code has weight at least
$p^2+1-t$. Thus $C_{g,r,w}$ has $(t+1,p^2+1-t)$ locality,
and up to $p^2-t$ \textit{erasures} \rev{(missing coordinate values)} within $\PP_C$ can be recovered
using its remaining coordinates.
Each coordinate has at least
\[
 \left\lfloor\frac{p^2}{t+1}\right\rfloor
\]
pairwise disjoint recovery sets, each of size $t+1$.
In particular, the families in Proposition~\ref{prop:explicit} are
asymptotically good with locality bounded independently of $r$.
\end{corollary}
\begin{proof}
Each $C$ is a projective line over $\overline{k}$,
and $\omega_{\overline{k}}^w|_C$ has degree $t$.
The space of sections of this line bundle has dimension $t+1$,
which bounds the dimension of the restricted code after extension
of scalars, and hence also over $k$.
A nonzero restriction of $f_{\overline{k}}$, with
$f\in M_{g,r,w}$, therefore
has at most $t$ zeros among the $p^2+1$ evaluation points, proving
the bound for the weight and the assertion about erasures.

Every evaluation point $P$ lies on such a curve. Since
$t+1\le p^2$, choose $t+1$ distinct points of $\PP_C\setminus\{P\}$.
If $f\in M_{g,r,w}$ evaluates to zero at all of them, then
$f_{\overline{k}}|_C=0$, and hence $f(P)=0$. For the
$k$-linear code, this kernel inclusion says that the coordinate at
$P$ is a $k$-linear combination of the chosen coordinates. 
This holds for every subset of size $t+1$ in
$\PP_C\setminus\{P\}$. Taking pairwise disjoint such subsets
gives the stated number of recovery sets.
\end{proof}

\section{Further questions}\label{sec:directions}

Theorem~\ref{thm:general} leaves several quantitative questions.
An effective estimate for the twists $b_i$ in
\eqref{eq:auxiliary-classes} would give a numerical characteristic
threshold. Proposition~\ref{prop:weighted-incidence} allows
different constants at successive dimensions. Computing the
degrees of the flag strata and using their individual incidence
numbers in place of $g^2$ may therefore improve the distance bound.

The range $0\le w\le p+1$ in Theorem~\ref{thm:mainideal}
comes from the restriction to superspecial points on a
supersingular curve in genus two. For $g=2,3$, vanishing on the
entire supersingular locus implies divisibility by $\Ha$ in the
larger range $w<(p^2-1)/2$.
It would be interesting to determine the corresponding range
in higher genus. The cycle classes in \cite{vdGH} and the
geometric descriptions \rev{after perfection} in \cite{ShimadaTakamatsu} may be useful
for this question, but a cycle class alone does not bound
containment in a divisor when the locus has higher codimension.
For example, in genus four the locus of $p$-rank zero has
dimension six, whereas the supersingular locus has dimension four.
An extension of the argument in genus three would require control
of restrictions to suitable intermediate strata, including the
degrees and multiplicities of the morphisms used to describe them.

The estimate for dimensions in Section~\ref{sec:growth} includes
the smallest weight $g+1$ of a scalar holomorphic discrete series.
Treating weights $w\le g$ requires a different argument.
\rev{Likewise, the distance estimate in genus three uses an additional
positivity condition in the proof of Proposition~\ref{prop:distance3};} extending
the range of weights in Theorem~\ref{thm:mainideal} alone would
not extend that distance bound.

Finally, the proof using flags counts auxiliary geometric points.
It does not by itself produce the recovery sets of
\rev{Section~\ref{subsec:local-recovery}}. A construction of such sets in
arbitrary genus, and algorithms for global decoding of these
codes, would require additional arguments.
Computing generators efficiently is a further question.

\appendix
\section{The supersingular locus in genus three}\label{sec:genus3}
\subsection{Strata and cycle classes}
In this appendix $g=3$ and $r\ge1$ is fixed. We prove the
distance estimate used in Proposition~\ref{prop:explicit},
and the stronger assertion about vanishing on the supersingular
locus announced after Theorem~\ref{thm:mainideal}.
\rev{To pass from the geometric Hasse divisor
$H_{\overline{k}}\subset(\overline X_{3,r})_{\overline{k}}$
to the supersingular locus, we use the intermediate subvarieties
defined by the $p$-rank.}

For $i=1,2,3$, let $Z_i\subset\overline X_{3,r}$ be the reduced
closure over $k$ of the locus in $X_{3,r}$ where the $p$-rank
is at most $3-i$.
Each $Z_i$ has dimension $6-i$, and after base change we have
\[
 H_{\overline{k}}=(Z_1)_{\overline{k}}
 \supset(Z_2)_{\overline{k}}\supset(Z_3)_{\overline{k}}
 \supset(S_{3,r})_{\overline{k}}.
\]

For $i=1,2,3$, let $U_i\subset X_{3,r}$ be the reduced locally
closed subscheme where the $p$-rank is exactly $3-i$.
The inverse image of $(U_i)_{\overline{k}}$ under
$\pi:\A_{3,N\ell^r}\to(X_{3,r})_{\overline{k}}$
is precisely the locus of $p$-rank $3-i$ in $\A_{3,N\ell^r}$.
Each of these three loci is irreducible by
\cite[proof of Proposition~3.1, Theorem~3.2]{YuMono}.
Since $\pi$ is surjective, $(U_i)_{\overline{k}}$ is irreducible.
The closure relation for the $p$-rank strata identifies
the closure of $(U_i)_{\overline{k}}$ in $(X_{3,r})_{\overline{k}}$
with the locus of $p$-rank at most $3-i$;
see \cite[Section~2]{vdGH}.
Since $k$ is perfect, $(Z_i)_{\overline{k}}$ is therefore the
reduced closure of $(U_i)_{\overline{k}}$ in
$(\overline X_{3,r})_{\overline{k}}$.
Thus $(Z_i)_{\overline{k}}$ is integral, and $Z_i$ is geometrically integral.
In particular, $\dim(Z_3)_{\overline{k}}=3$, while every
irreducible component of $(S_{3,r})_{\overline{k}}$ has dimension two.
Since $(S_{3,r})_{\overline{k}}\subset(Z_3)_{\overline{k}}$,
each such component is a prime divisor on $(Z_3)_{\overline{k}}$.

For the Chern classes $\lambda_{i,\overline{k}}$ of the extended Hodge bundle $\mathbb E_{\overline{k}}$
on $(\overline X_{3,r})_{\overline{k}}$, with $\lambda_{1,\overline{k}}=\lambda_{\overline{k}}$, we use the following relations:
\begin{align}
 [(Z_i)_{\overline{k}}]&=\prod_{j=1}^{i}(p^j-1)\lambda_{i,\overline{k}},
 \label{eq:prank}\\
 \lambda_{1,\overline{k}}^2&=2\lambda_{2,\overline{k}},\qquad
 \lambda_{2,\overline{k}}^2=2\lambda_{1,\overline{k}}\lambda_{3,\overline{k}},\qquad
 \lambda_{1,\overline{k}}^4=8\lambda_{1,\overline{k}}\lambda_{3,\overline{k}},
 \label{eq:taut}\\
 [(S_{3,r})_{\overline{k}}]&=(p-1)^2(p^3-1)(p^4-1)\lambda_{1,\overline{k}}\lambda_{3,\overline{k}}
      =(p-1)(p^2+1)\lambda_{\overline{k}}[(Z_3)_{\overline{k}}].
 \label{eq:class3}
\end{align}
These identities hold in the Chow ring of
$(\overline X_{3,r})_{\overline{k}}$ with rational coefficients.
Equation~\eqref{eq:prank} is the formula for the classes of
the $p$-rank loci in \cite[Section~2, formula~(1)]{vdGH}.
The relations in \eqref{eq:taut} follow from the tautological
relation $c(\mathbb E_{\overline{k}})c(\mathbb E_{\overline{k}}^\vee)=1$;
see \cite[Section~2]{vdGH}.
The first equality in \eqref{eq:class3} is
\cite[Theorem~8.1]{vdGH}, and the second follows from
\eqref{eq:prank} with $i=3$.

\subsection{Ruled models and incidence relations}

To relate restriction to the surfaces in $(S_{3,r})_{\overline{k}}$ to evaluation
at superspecial points, we use curves defined by the $a$-number.
For an abelian variety $A$ over $\overline{k}$, its \textit{$a$-number} is defined by
$a(A)=\dim_{\overline{k}}\operatorname{Hom}(\alpha_p,A[p])$,
where $\alpha_p$ is the kernel of Frobenius on the additive group.
Let $\YY$ be the set of irreducible components of $(S_{3,r})_{\overline{k}}$, and let
$\CC$ be the set of one-dimensional irreducible components
of the reduced locus in $(S_{3,r})_{\overline{k}}$ defined by $a(A)\ge2$.
These curves are the irreducible components of the closure
of the Ekedahl--Oort stratum \rev{indexed by the Young diagram
$[3,2]$} in the notation
of \cite[Section~9]{vdGH}.
The following proposition records their degrees with respect to
the Hodge line bundle and their incidence relations.
\begin{proposition}\label{prop:geometry3}
\begin{enumerate}
    \item  At every level $r$, each $Y\in\YY$ contains exactly $p^3+1$ curves
in $\CC$, and $\lambda_{\overline{k}}^2[Y]=2(p-1)^2(p+1)$.

\item Each $C\in\CC$ is isomorphic to $\mathbb P^1_{\overline{k}}$,
satisfies $\lambda_{\overline{k}}[C]=p-1$, and contains exactly
$p^2+1$ superspecial points.

\item Each curve lies on $p^2+1$ surfaces, and each superspecial point
lies on $p^3+1$ curves.
\end{enumerate}
\end{proposition}
\begin{proof}
We first construct a ruled model for each $Y\in\YY$, and
then use this model to compute the degrees and incidences.
The flag construction in
\cite[Definition~4.1 and Section~7]{vdGH}
starts with a polarized superspecial threefold $(A_0,\eta_0)$ over $\overline{k}$
\rev{such that
\[
 \ker\eta_0=A_0[F^2]=A_0[p],\qquad \deg\eta_0=p^6.
\]
Here $F^2$ is the composite of two successive relative Frobenius
morphisms. The construction}
parametrizes chains of isogenies
$A_0\times_{\overline{k}} B\to\mathscr A_1\to\mathscr A$.
Here $\mathscr A_1$ and $\mathscr A$ are abelian schemes over
\[
 B=\mathbb P\bigl(\mathcal O_T(1)\oplus\mathcal O_T(-1)\bigr),
 \qquad
 T:\ x^{p+1}+y^{p+1}+z^{p+1}=0,
\]
where $T$ is the \textit{Hermitian curve}, viewed over $\overline{k}$.
\rev{Here $\mathcal O_T(1)=\mathcal O_{\mathbb P^2}(1)|_T$.}
The terminal abelian scheme $\mathscr A$ carries a principal
polarization $\eta$.
If $q:A_0\times_{\overline{k}} B\to\mathscr A$ is the composite isogeny,
then $q$ has $p$-power degree and
$q^*\eta=(\eta_0)_B$.

We now add the level structures defining $(X_{3,r})_{\overline{k}}$.
For $m=N$ and $m=\ell^r$, the isogeny $q$ induces an isomorphism
$q_m:A_0[m]\times_{\overline{k}} B\xrightarrow{\sim}\mathscr A[m]$,
because $p\nmid N\ell$.
The equality $q^*\eta=(\eta_0)_B$ makes these isomorphisms
compatible with the Weil pairings defined by $\eta_0$ and $\eta$.
Consequently, a full symplectic level $N$ structure on $A_0[N]$
and a Lagrangian direct summand of $A_0[\ell^r]$ determine
unique corresponding level structures on $\mathscr A$.
Each such choice therefore gives a morphism
$B\to(X_{3,r})_{\overline{k}}$.

To determine its degree, put $G=\operatorname{Aut}(A_0,\eta_0)$,
and let $\mathcal L$ be the finite set of these choices of
level structures on $A_0$.
By \cite[Section~7 and proof of Theorem 8.1]{vdGH}, the morphism from the quotient stack
$[B/G]$ to the corresponding supersingular component of the
moduli stack has degree one.
\rev{Base change by the finite \'etale level cover gives the source
$(B\times\mathcal L)/G$.}
The action of $G$ on $\mathcal L$ is free since $N\geq 3$.
Choosing a representative of each $G$-orbit in $\mathcal L$
thus identifies this quotient with a disjoint union of
copies of $B$.
The fact that this degree is one therefore gives, for every $Y\in\YY$,
a proper morphism $\nu:B\to Y\subset(X_{3,r})_{\overline{k}}$
that is generically finite of degree one. In particular,
\[
 [\overline{k}(B):\overline{k}(Y)]=1,\qquad \nu_*[B]=[Y].
\]
This is the full degree of the extension of function fields,
including its inseparable degree.
The morphism contracts the negative section of the ruled
surface, as in \cite[Section~7]{vdGH}.

We next compute the degree of $Y$ with respect to $\omega_{\overline{k}}$.
Let $\ell_0,\ell_1$ be the classes defined in
\cite[Section~8]{vdGH}.
Here $\ell_1$ is pulled back from the hyperplane class on $T$,
and $\ell_0$ has degree one on a fiber of $B\to T$.
The calculation of the Hodge bundle therein gives
$\nu^*\lambda_{\overline{k}}=(p-1)(\ell_0+\ell_1)$,
$\ell_0^2=\ell_1^2=0$, and $\deg(\ell_0\ell_1)=p+1$.
Using $\nu_*[B]=[Y]$ and the projection formula, we obtain
\[
 \lambda_{\overline{k}}^2[Y]
 =(\nu^*\lambda_{\overline{k}})^2[B]
 =(p-1)^2(\ell_0+\ell_1)^2[B]
 =2(p-1)^2(p+1).
\]

The curves in $\CC$ contained in $Y$ are described in
\cite[Section~9]{vdGH}.
They are the distinct images of the fibers of $B\to T$
over the points of $T$ with projective coordinates in $\F_{p^2}$
for the displayed Hermitian equation.
Each image is isomorphic to $\mathbb P^1_{\overline{k}}$,
has degree $p-1$ with respect to $\omega_{\overline{k}}$,
and contains exactly $p^2+1$ superspecial points.
There are $p^3+1$ such points of $T$, and hence exactly $p^3+1$
such curves in $Y$.

Finally, the numbers of incident flags computed in
\cite[Section~9]{vdGH} give $p^2+1$ surfaces through a curve
and $p^3+1$ curves through a superspecial point.
These counts apply with our level structures as follows.
Fix the underlying polarized abelian variety at a geometric
generic point of a curve, or at a superspecial point,
together with its level structures.
Along each incident flag, these level structures pull back
uniquely, since every isogeny in the flag induces isomorphisms
on the $N$-torsion and the $\ell^r$-torsion.
Conversely, forgetting the transported level structures
recovers the original flag.
The chosen full level $N$ structure eliminates automorphisms
of the fixed polarized object, so the same flag counts give
the stated numbers of incident components at level $r$.
All superspecial points counted here are points of
$(X_{3,r})_{\overline{k}}$; the contracted section contributes
only its image point.
\end{proof}

For $Y\in\YY$, write $\nu_Y:B_Y\to Y$ for the ruled model above,
$E_Y$ for its contracted section, and $P_Y=\nu_Y(E_Y)$.
By \cite[Proposition~8.2]{vdGH}, $E_Y^2=-2(p+1)$.
Each of the $p^3+1$ specified fibers maps isomorphically to its
image in $\CC$: the degrees of $\nu_Y^*\omega_{\overline{k}}$
on the fiber and of $\omega_{\overline{k}}$ on its image
are both $p-1$, so the map between the two projective
lines has degree one.

\begin{lemma}\label{lem:contraction3}
The map $Y\mapsto P_Y$ is a bijection $\YY\to\PP_{3,r}$.
For $Y\in\YY$ and $C\in\CC$, one has $C\subset Y$ if and only if
$P_Y\in C$.
\end{lemma}
\begin{proof}
 The contracted point is superspecial by
the flag description in \cite[Sections~7 and~9]{vdGH}.
The two incidence counts in Proposition~\ref{prop:geometry3} give
$|\YY|(p^3+1)=|\CC|(p^2+1)=|\PP_{3,r}|(p^3+1)$, so $|\YY|=|\PP_{3,r}|$.
We prove that every $P\in\PP_{3,r}$ is a contracted point.

Represent the geometric image of $P$ by a principally polarized
superspecial threefold $(A,\lambda_A)$ over $\overline{k}$ with its level structures. Set
$(A',\lambda_{A'})=(A,\lambda_A)^{(p^{-1})}$, where $(p^{-1})$
denotes pullback by the inverse of Frobenius on $\overline{k}$, and let
$F:A'\to A$ be relative Frobenius. Then
$F^*\lambda_A=p\lambda_{A'}$, whose kernel is
$A'[p]=A'[F^2]$. Here $F^2$ denotes the composite of the two
successive relative Frobenius morphisms from $A'$ to $A'^{(p^2)}$.
The latter equality follows from a superspecial
model with $p^2$-Frobenius $[-p]$, as in
\cite[Proposition~5.1 and Definition~5.2]{YuArithmetic}.
Pull back the level structures along $F$, which is an isomorphism
on prime-to-$p$ torsion. The polarized source with these compatible
levels indexes a supersingular component by
\cite[Sections~3 and~7]{vdGH} and Proposition~\ref{prop:geometry3}.
In its flag family, the intermediate quotients correspond to
subgroups $K\simeq\alpha_p^2$ of $A'[F]$ on the Hermitian curve.
The quotients
$A'\to A'/K\to A'/A'[F]=A$
define a section of the ruled model with constant polarized target
and constant level structure. It is the contracted section.
Its image is $P$, proving surjectivity and hence bijectivity.

If $C\subset Y$, its corresponding fiber meets $E_Y$, so $P_Y\in C$.
The $p^2+1$ surfaces containing a fixed $C$ therefore give
$p^2+1$ distinct superspecial points on $C$. These are all such
points by Proposition~\ref{prop:geometry3}, proving the converse.
\end{proof}

\subsection{Vanishing on the supersingular locus}

We first establish nonvanishing on the integral threefold
$(Z_3)_{\overline{k}}$.

\begin{lemma}\label{lem:restriction3}
Let $g=3$, $0\le w<(p^2-1)/2$, and $f\in M_{3,r,w}$.
If $f$ is not divisible by $\Ha$, then
$f_{\overline{k}}|_{(Z_i)_{\overline{k}}}\ne0$ for $i=1,2,3$.
\end{lemma}

\begin{proof}
Since $Z_1=H$, equation~\eqref{eq:hasse-restriction} gives
$f|_{Z_1}\ne0$ over $k$.
By faithful flatness, $f_{\overline{k}}|_{(Z_1)_{\overline{k}}}\ne0$.
The chain
$(Z_1)_{\overline{k}}\supset(Z_2)_{\overline{k}}\supset(Z_3)_{\overline{k}}$
consists of integral subvarieties of dimensions five, four, and three.
Put $v_i=\lambda_{\overline{k}}^{6-i}[(Z_i)_{\overline{k}}]$
for $i=1,2,3$.

The tautological identities \eqref{eq:taut} imply
$\lambda_{\overline{k}}^4\lambda_{2,\overline{k}}=\lambda_{\overline{k}}^6/2$ and
$\lambda_{\overline{k}}^3\lambda_{3,\overline{k}}=\lambda_{\overline{k}}^6/8$.
Combining these with the cycle classes \eqref{eq:prank} gives
\begin{equation}\label{eq:degrees3}
 \begin{aligned}
  v_1&=(p-1)\lambda_{\overline{k}}^6,\\
  v_2&=\frac{(p-1)(p^2-1)}2\lambda_{\overline{k}}^6,\\
  v_3&=\frac{(p-1)(p^2-1)(p^3-1)}8\lambda_{\overline{k}}^6.
 \end{aligned}
\end{equation}
All three degrees are positive, since $\lambda_{\overline{k}}^6>0$
by Section~\ref{subsec:hasse} and invariance of the top intersection
degree under extension of the ground field.
Their successive ratios are $v_2/v_1=(p^2-1)/2$ and
$v_3/v_2=(p^3-1)/4$. Since $p\ge5$, the weight hypothesis gives
$w<(p^2-1)/2<(p^3-1)/4$.
If the restriction first vanished on $Z_i$, Lemma~\ref{lem:intersection}
on $Z_{i-1}$ would give $v_i\le wv_{i-1}$, contradicting
these inequalities. Thus both restrictions are nonzero.
\end{proof}

\begin{proposition}\label{prop:ideal3}
For $g=3$ and $0\le w<(p^2-1)/2$, the following equalities hold
as subspaces of $M_{3,r,w}$ over $k$:
\[
 \begin{aligned}
 H^0(\overline X_{3,r},\I_{S_{3,r}}\otimes\omega^w)
 &=H^0(\overline X_{3,r},\I_{Z_3}\otimes\omega^w)\\
 &=H^0(\overline X_{3,r},\I_{Z_2}\otimes\omega^w)\\
 &=\Ha M_{3,r,w-(p-1)}.
 \end{aligned}
\]
\end{proposition}

\begin{proof}
Every multiple of $\Ha$ vanishes on all the indicated loci,
since $S_{3,r}\subset Z_3\subset Z_2\subset H$ over $k$.
Conversely, suppose that
$f\in H^0(\overline X_{3,r},\I_{S_{3,r}}\otimes\omega^w)$
is not divisible by $\Ha$.
Lemma~\ref{lem:restriction3} gives
$f_{\overline{k}}|_{(Z_3)_{\overline{k}}}\ne0$.
Each irreducible component $Y\in\YY$ of $(S_{3,r})_{\overline{k}}$
is an integral surface, hence an integral divisor on $(Z_3)_{\overline{k}}$,
and $f_{\overline{k}}|_{(Z_3)_{\overline{k}}}$ vanishes on every such $Y$.

Applying Lemma~\ref{lem:intersection} on $(Z_3)_{\overline{k}}$ and then
using \eqref{eq:class3}, we obtain
\[
 \begin{aligned}
  w\lambda_{\overline{k}}^3[(Z_3)_{\overline{k}}]
  &\ge\sum_{Y\in\YY}\lambda_{\overline{k}}^2[Y]
   =\lambda_{\overline{k}}^2[(S_{3,r})_{\overline{k}}]\\
  &=(p-1)(p^2+1)\lambda_{\overline{k}}^3[(Z_3)_{\overline{k}}].
 \end{aligned}
\]
Here the sum equals $\lambda_{\overline{k}}^2[(S_{3,r})_{\overline{k}}]$
because $(S_{3,r})_{\overline{k}}$ is reduced.
By \eqref{eq:degrees3}, $\lambda_{\overline{k}}^3[(Z_3)_{\overline{k}}]>0$.
Dividing by this degree yields
$w\ge(p-1)(p^2+1)>(p^2-1)/2$,
contrary to the weight hypothesis.
This proves
\[
 H^0(\overline X_{3,r},\I_{S_{3,r}}\otimes\omega^w)
 =\Ha M_{3,r,w-(p-1)}.
\]

Finally, $S_{3,r}\subset Z_3\subset Z_2\subset H$ over $k$ implies
\[
 \begin{aligned}
 \Ha M_{3,r,w-(p-1)}
 &\subset H^0(\overline X_{3,r},\I_{Z_2}\otimes\omega^w)\\
 &\subset H^0(\overline X_{3,r},\I_{Z_3}\otimes\omega^w)\\
 &\subset H^0(\overline X_{3,r},\I_{S_{3,r}}\otimes\omega^w).
 \end{aligned}
\]
All these inclusions are therefore equalities.
\end{proof}

\subsection{The distance bound}

\begin{proposition}\label{prop:distance3}
For $g=3$ and $0\le w\le(p-1)/2$, every nonzero $f\in M_{3,r,w}$
satisfies
\[
 \frac{\wt(\ev(f))}{n_{3,r}}
 \ge\left(1-\frac{w}{(p-1)(p^2+1)}\right)
 \frac{(p^2+1-w(p-1))(p^2-p-1)}
 {(p^2+1)(p^2-p-1)+2w(p-1)}>0.
\]
In particular, evaluation is injective.
\end{proposition}

\begin{proof}
The case $w=0$ follows from $M_{3,r,0}=k$. Assume $w\ge1$ and fix
$0\ne f\in M_{3,r,w}$. Since $w<p-1$, the absence of modular forms
of negative weight and Lemma~\ref{lem:restriction3} give $f_{\overline{k}}|_{(Z_3)_{\overline{k}}}\ne0$.
Let $u_2$ and $u_1$ be the proportions of members of $\YY$ and
$\CC$, respectively, on which $f_{\overline{k}}$ is not identically zero, and
put $\delta=\wt(\ev(f))/n_{3,r}$. By Lemma~\ref{lem:contraction3},
$|\YY|=n_{3,r}$.

Since the surfaces in $\YY$ have the same degree with respect to
$\omega_{\overline{k}}$, Lemma~\ref{lem:intersection} on
$(Z_3)_{\overline{k}}$ gives
\[
 (1-u_2)\lambda_{\overline{k}}^2[(S_{3,r})_{\overline{k}}]
 \le w\lambda_{\overline{k}}^3[(Z_3)_{\overline{k}}].
\]
Using \eqref{eq:class3}, we obtain $u_2\ge1-\alpha$, where
$\alpha=w/((p-1)(p^2+1))$.
Set
$t=w(p-1)$, $b=p^2+1$, $h=p^3+1$, $e=2(p+1)$, and
$\gamma=t/b$.
\rev{The weight hypothesis gives $t\le(p-1)^2/2<b$, hence
$0\le\gamma<1$.}
Then $|\CC|b=n_{3,r}h$.
On each curve with nonzero restriction, the zero divisor has
degree $t$. Counting incidences with nonzero evaluation points
therefore gives
\begin{equation}\label{eq:curve-support3}
 \delta\ge(1-\gamma)u_1.
\end{equation}

For a surface $Y$ with $f_{\overline{k}}|_Y\ne0$,
let $z_Y$ be the number of curves in $Y$ on which
$f_{\overline{k}}$ vanishes identically.
If $f(P_Y)\ne0$, then $z_Y=0$ by
Lemma~\ref{lem:contraction3}.
Suppose now that $f(P_Y)=0$.
Since $\nu_Y$ contracts $E_Y$ to $P_Y$, the zero
divisor of $\nu_Y^*f_{\overline{k}}$ contains $E_Y$.
We may therefore write
$\operatorname{div}(\nu_Y^*f_{\overline{k}})=a_YE_Y+R_Y$
with $a_Y\ge1$,
where $R_Y$ is effective and does not contain $E_Y$.
Let $F$ be a general fiber of $B_Y\to T$.
Since $E_Y\cdot F=1$ and $R_Y\cdot F\ge0$,
\[
t=(a_YE_Y+R_Y)\cdot F
=a_Y+R_Y\cdot F\ge a_Y.
\]
Moreover, $\nu_Y^*\omega_{\overline{k}}^w$
has degree zero on $E_Y$, so
\[
0=(a_YE_Y+R_Y)\cdot E_Y
=-ea_Y+R_Y\cdot E_Y.
\]
Thus $R_Y\cdot E_Y=ea_Y$.
Each of the $z_Y$ curves on which $f_{\overline{k}}$ vanishes identically corresponds
to a fiber contained in $R_Y$, and each such
fiber meets $E_Y$ once.
Since $R_Y$ is effective and does not contain
$E_Y$, we obtain
$z_Y\le R_Y\cdot E_Y=ea_Y$.

For each $C\subset Y$ with
$f_{\overline{k}}|_C\ne0$, the corresponding
fiber maps isomorphically onto $C$ and meets
$E_Y$ once.
\rev{On the smooth surface $B_Y$, choose a local equation $u$ for
$E_Y$ at this intersection point and a local trivialization of
$\nu_Y^*\omega_{\overline{k}}^w$. Then
$\nu_Y^*f_{\overline{k}}=u^{a_Y}v$ for a regular function $v$.
Since $u$ restricts to a local parameter on the fiber, the nonzero
restriction to that fiber has vanishing order at least $a_Y$.}
Hence
\begin{equation}\label{eq:orders3}
\operatorname{ord}_{P_Y}(f_{\overline{k}}|_C)
\ge a_Y.
\end{equation}

We next combine the bounds on $z_Y$ and the vanishing orders
at $P_Y$ to relate $u_2$, $u_1$, and $\delta$.
Set $\xi=h/e-t-1$.
The weight assumption gives
$\xi\ge(p-2)/2>0$.
For $1\le a\le t$ and $0\le z\le ea$, we have
\[
a(h-z)-et-\xi z
=e(a-1)(t-a)+(ea-z)(a+\xi)\ge0.
\]

Let $\YY_0$ be the set of surfaces $Y\in\YY$
such that $f_{\overline{k}}|_Y\ne0$ and
$f(P_Y)=0$.
Since $Y\mapsto P_Y$ is a bijection and
$f(P_Y)\ne0$ implies $f_{\overline{k}}|_Y\ne0$,
we have
$|\YY_0|=n_{3,r}(u_2-\delta)$.

Applying the preceding inequality to each
$Y\in\YY_0$ and summing gives
\[
et|\YY_0|+\xi\sum_{Y\in\YY_0}z_Y
\le
\sum_{Y\in\YY_0}a_Y(h-z_Y).
\]

We first count the curves on which $f_{\overline{k}}$ vanishes
identically inside the surfaces on which it does not vanish
identically. We have $z_Y=0$ whenever $f(P_Y)\ne0$.
Moreover, every curve on which $f_{\overline{k}}$ is nonzero
lies on exactly $b$ such surfaces. Hence
\[
\sum_{Y\in\YY_0}z_Y
=n_{3,r}hu_2-|\CC|bu_1
=n_{3,r}h(u_2-u_1).
\]

Next, fix a curve $C\in\CC$ with
$f_{\overline{k}}|_C\ne0$.
By Lemma~\ref{lem:contraction3}, the points
$P_Y$ for surfaces $Y$ containing $C$
are precisely the $b$ superspecial points
on $C$.
Using \eqref{eq:orders3}, we obtain
\[
\sum_{\substack{Y\in\YY_0\\ C\subset Y}}a_Y
\le\sum_{P\in C\cap\PP_{3,r}}
\operatorname{ord}_P(f_{\overline{k}}|_C)
\le t,
\]
where the last inequality follows from
$\deg(\omega_{\overline{k}}^w|_C)=t$.

Summing over all $C\in\CC$ with $f_{\overline{k}}|_C\ne0$ gives
\[
\sum_{Y\in\YY_0}a_Y(h-z_Y)
\le t|\CC|u_1.
\]
Consequently, using
$|\CC|b=n_{3,r}h$ and $\gamma=t/b$,
we obtain
$et(u_2-\delta)+\xi h(u_2-u_1)\le\gamma h u_1$.
Use $u_1\le\delta/(1-\gamma)$ from
\eqref{eq:curve-support3}. Since $\xi+\gamma>0$, rearrangement gives
\[
 \delta\ge
 \frac{(1-\gamma)(et+\xi h)}{et(1-\gamma)+h(\xi+\gamma)}u_2
 =(1-\gamma)\frac{h-e}{h-e+e\gamma}u_2.
\]
Here $h-et>0$, and the equality follows by substituting
$\xi=h/e-t-1$ and cancelling $h-et$.
Finally use $u_2\ge1-\alpha$ and
$h-e=(p+1)(p^2-p-1)$ to obtain the stated bound.
Its factors are positive in the given weight range, so every
nonzero $f$ has nonzero evaluation.
\end{proof}

Counting only the curves on which $f_{\overline{k}}$ vanishes
identically inside each surface gives the product bound
$(1-\alpha)(1-2w(p^2-1)/(p^3+1))(1-\gamma)$.
Proposition~\ref{prop:distance3} strengthens it by counting
vanishing orders along the contracted sections.

\section*{Acknowledgments}
The author thanks Takuya Yamauchi for his helpful comments.
The author used ChatGPT (OpenAI) in the preparation of this manuscript.
The author takes full responsibility for its content, including the
correctness of the mathematical statements.

\bibliographystyle{amsalpha}
\bibliography{main}

@article{TVZ,
 author={Tsfasman, M. A. and Vl{\u a}du{\c t}, S. G. and Zink, T.},
 title={Modular curves, {Shimura} curves, and {Goppa} codes, better than {Varshamov--Gilbert} bound},
 journal={Mathematische Nachrichten}, volume={109}, year={1982}, pages={21--28},
 doi={10.1002/mana.19821090103}}

@incollection{ChaiOort,
 author={Chai, Ching-Li and Oort, Frans},
 title={Moduli of abelian varieties and $p$-divisible groups},
 booktitle={Arithmetic Geometry},
 editor={Darmon, Henri and Ellwood, David and Hassett, Brendan and Tschinkel, Yuri},
 series={Clay Mathematics Proceedings}, volume={8},
 publisher={American Mathematical Society}, address={Providence, RI},
 year={2009}, pages={441--536},
 url={https://www.claymath.org/wp-content/uploads/2022/03/cmip08c.pdf}}

@book{FC,
 author={Faltings, Gerd and Chai, Ching-Li}, title={Degeneration of Abelian Varieties},
 series={Ergebnisse der Mathematik und ihrer Grenzgebiete}, volume={22},
 publisher={Springer-Verlag}, address={Berlin}, year={1990}, doi={10.1007/978-3-662-02632-8}}

@article{KO,
 author={Katsura, Toshiyuki and Oort, Frans}, title={Families of supersingular abelian surfaces},
 journal={Compositio Mathematica}, volume={62}, year={1987}, pages={107--167},
 url={https://www.numdam.org/item/CM_1987__62_2_107_0/}}

@article{YuMass,
 author={Yu, Chia-Fu}, title={The supersingular loci and mass formulas on {Siegel} modular varieties},
 journal={Documenta Mathematica}, volume={11}, year={2006}, pages={449--468},
 url={https://arxiv.org/abs/math/0608458}}

@article{YuMono,
 author={Yu, Chia-Fu}, title={Irreducibility and $p$-adic monodromies on the {Siegel} moduli spaces},
 journal={Advances in Mathematics}, volume={218}, year={2008}, pages={1253--1285},
 url={https://webdoc.sub.gwdg.de/ebook/serien/e/mpi_mathematik/2007/3.pdf}}

@incollection{vdG,
 author={van der Geer, Gerard}, title={Cycles on the moduli space of abelian varieties},
 booktitle={Moduli of Curves and Abelian Varieties}, series={Aspects of Mathematics},
 volume={E33}, publisher={Vieweg}, year={1999}, pages={65--89},
 url={https://arxiv.org/abs/alg-geom/9605011}}

@article{vdGH,
 author={van der Geer, Gerard and Harashita, Shushi},
 title={The cycle class of the supersingular locus of principally polarized abelian varieties},
 journal={Selecta Mathematica. New Series}, volume={31}, year={2025}, pages={Paper No.~95},
 doi={10.1007/s00029-025-01095-w}, url={https://arxiv.org/abs/2307.14393}}

@article{FinisLapid,
 author={Finis, Tobias and Lapid, Erez},
 title={On the analytic properties of intertwining operators {II}: local degree bounds and limit multiplicities},
 journal={Israel Journal of Mathematics}, volume={230}, year={2019}, pages={771--793},
 url={https://arxiv.org/abs/1705.08191}}

@article{Wakatsuki,
 author={Wakatsuki, Satoshi}, title={The dimensions of spaces of {Siegel} cusp forms of general degree},
 journal={Advances in Mathematics}, volume={340}, year={2018}, pages={1012--1066},
 url={https://arxiv.org/abs/1602.05676}}

@article{Hansen,
 author={Hansen, S{\o}ren Have}, title={Error-correcting codes from higher-dimensional varieties},
 journal={Finite Fields and Their Applications}, volume={7}, year={2001}, pages={530--552},
 doi={10.1006/ffta.2001.0313}}

@incollection{CLP,
 author={Couvreur, Alain and Lebacque, Philippe and Perret, Marc},
 title={Toward good families of codes from towers of surfaces},
 booktitle={Arithmetic, Geometry, Cryptography and Coding Theory}, series={Contemporary Mathematics},
 volume={770}, publisher={American Mathematical Society}, year={2021}, pages={59--101},
 url={https://arxiv.org/abs/2002.02220}}

@unpublished{BlacheHallouin,
 author={Blache, R{\'e}gis and Hallouin, Emmanuel},
 title={Ruled surfaces over finite fields, and some codes over them},
 note={Preprint, arXiv:2509.18698}, year={2025},
 url={https://arxiv.org/abs/2509.18698}}

@article{Ghitza,
 author={Ghitza, Alexandru}, title={Hecke eigenvalues of {Siegel} modular forms (mod $p$) and of algebraic modular forms},
 journal={Journal of Number Theory}, volume={106}, year={2004}, pages={345--384},
 url={https://arxiv.org/abs/math/0309006}}

@article{Yamauchi,
 author={Yamauchi, Takuya}, title={The weight reduction of mod $p$ {Siegel} modular forms for {$GSp_4$} and theta operators},
 journal={Mathematische Zeitschrift}, volume={303}, year={2023}, pages={Paper No.~10},
 doi={10.1007/s00209-022-03153-x}, url={https://arxiv.org/abs/1410.7894}}

@unpublished{ShimadaTakamatsu,
 author={Shimada, Ryosuke and Takamatsu, Teppei},
 title={On the supersingular locus of the {Siegel} modular variety of genus 3 or 4},
 note={To appear in Annales de l'Institut Fourier; arXiv:2403.19505v2},
 year={2025}, url={https://arxiv.org/abs/2403.19505}}

@article{YuArithmetic,
 author={Yu, Chia-Fu}, title={On arithmetic of the superspecial locus},
 journal={Indiana University Mathematics Journal}, volume={67}, number={4},
 year={2018}, pages={1341--1382}, doi={10.1512/iumj.2018.67.7436},
 url={https://arxiv.org/abs/1210.1120}}

@article{FinisLapidApprox,
  author  = {Finis, Tobias and Lapid, Erez},
  title   = {An approximation principle for congruence subgroups
             {II}: application to the limit multiplicity problem},
  journal = {Mathematische Zeitschrift},
  volume  = {289},
  number  = {3-4},
  pages   = {1357--1380},
  year    = {2018},
  doi     = {10.1007/s00209-017-2002-0},
  url     = {https://arxiv.org/abs/1504.04795}
}

@article{BargTamoVladut,
  author  = {Barg, Alexander and Tamo, Itzhak
             and Vl{\u{a}}du{\c{t}}, Serge},
  title   = {Locally recoverable codes on algebraic curves},
  journal = {IEEE Transactions on Information Theory},
  volume  = {63},
  number  = {8},
  pages   = {4928--4939},
  year    = {2017},
  doi     = {10.1109/TIT.2017.2700859},
  url     = {https://arxiv.org/abs/1603.08876}
}

@incollection{BargEtAl,
  author    = {Barg, Alexander and Haymaker, Kathryn
               and Howe, Everett W. and Matthews, Gretchen L.
               and V{\'a}rilly-Alvarado, Anthony},
  title     = {Locally recoverable codes from algebraic curves
               and surfaces},
  booktitle = {Algebraic Geometry for Coding Theory and Cryptography},
  editor    = {Howe, Everett W. and Lauter, Kristin E.
               and Walker, Judy L.},
  publisher = {Springer},
  address   = {Cham},
  pages     = {95--127},
  year      = {2017},
  doi       = {10.1007/978-3-319-63931-4_4},
  url       = {https://arxiv.org/abs/1701.05212}
}

@article{CalegariGeraghty,
  author  = {Frank Calegari and David Geraghty},
  title   = {Minimal modularity lifting for nonregular
             symplectic representations},
  journal = {Duke Mathematical Journal},
  volume  = {169},
  number  = {5},
  pages   = {801--896},
  year    = {2020}
}

@incollection{EkedahlGeer,
  author    = {Ekedahl, Torsten and van der Geer, Gerard},
  title     = {Cycle classes of the {E-O} stratification on the
               moduli of abelian varieties},
  booktitle = {Algebra, Arithmetic, and Geometry:
               In Honor of {Yu. I. Manin}, Vol. I},
  editor    = {Tschinkel, Yuri and Zarhin, Yuri},
  series    = {Progress in Mathematics},
  volume    = {269},
  publisher = {Birkh{\"a}user Boston},
  year      = {2009},
  pages     = {567--636},
  doi       = {10.1007/978-0-8176-4745-2_13},
  url       = {https://van-der-geer.nl/strata_PiM269.pdf}
}

@unpublished{DAddezioVanHoften,
  author = {D'Addezio, Marco and van Hoften, Pol},
  title = {Hecke orbits on {Shimura} varieties of {Hodge} type},
  note = {Preprint, arXiv:2205.10344v2},
  year = {2025},
  url = {https://arxiv.org/abs/2205.10344v2}
}

@article{LanKoecher,
  author = {Lan, Kai-Wen},
  title = {Higher {Koecher}'s principle},
  journal = {Mathematical Research Letters},
  volume = {23},
  number = {1},
  year = {2016},
  pages = {163--199},
  doi = {10.4310/MRL.2016.v23.n1.a9},
  url = {https://www.kwlan.org/articles/Koecher.pdf}
}

@book{LanCompactifications,
  author = {Lan, Kai-Wen},
  title = {Arithmetic Compactifications of {PEL}-Type {Shimura} Varieties},
  series = {London Mathematical Society Monographs Series},
  volume = {36},
  publisher = {Princeton University Press},
  address = {Princeton, NJ},
  year = {2013},
  doi = {10.23943/princeton/9780691156545.001.0001},
  url = {https://www.kwlan.org/articles/cpt-PEL-type-thesis-revision.pdf}
}

@article{AtiyahSchmid,
  author = {Atiyah, Michael and Schmid, Wilfried},
  title = {A geometric construction of the discrete series for semisimple {Lie} groups},
  journal = {Invent. Math.},
  volume = {42},
  number = {1},
  pages = {1--62},
  year = {1977},
  doi = {10.1007/BF01389783},
  url = {https://www.math.utah.edu/~ptrapa/math-library/schmid/atiyah-schmid.pdf}
}
\end{document}